\documentclass{amsart}

\newtheorem{thm}{Theorem}[subsection]
\newtheorem{cor}[thm]{Corollary}
\newtheorem{lem}[thm]{Lemma}
\newtheorem{prop}[thm]{Proposition}
\theoremstyle{definition}
\newtheorem{defn}[thm]{Definition}
\newtheorem{example}[thm]{Example}

\theoremstyle{remark}
\newtheorem{rem}[thm]{Remark}
\newtheorem{Question}[thm]{Question}

\numberwithin{equation}{subsection}
\usepackage{amsmath,amssymb,graphicx,bbm,amsthm}
\begin{document}

\title[Alexander polynomials of ribbon knots and virtual knots]
{ Alexander polynomials of ribbon knots and virtual knots }

\author{Sheng Bai}
\address{  }
\email{ barries@163.com }



\begin{abstract}
We find that the Alexander polynomial of a ribbon knot in $ \mathbb{Z}HS^3 $ is determined by the intrinsic singularity information of its ribbon, and give explicit formulae to compute Alexander polynomial of the ribbon knot by that. We define half Alexander polynomial $  A_R (t) $, an invariant of oriented ribbons, and in fact the Alexander polynomial of the ribbon knot is $ A_R (t)  A_R (t^{-1}) $. We give two different simplified formulae for half Alexander polynomial. We characterize completely the polynomials arising as half Alexander polynomials of ribbons. 
The above study unexpectedly leads us to discover new formulae for the Alexander polynomial of general knots and virtual knots in terms of Gauss diagrams.
\end{abstract}

\date{2023 Jan.18}
\maketitle

\tableofcontents

\section{Introduction}
Alexander polynomial is one of the most famous and fundamental invariants of knots. There are many methods for computing Alexander polynomial. For example, construct a Seifert surface, choose a homological basis on it and calculate its Seifert matrix. The Alexander polynomial is also computed from a presentation of the knot group, the fundamental group of the knot complement. Thus one can compute the Alexander polynomial of a knot from its diagram, since there are many ways to get a presentation of the knot group from its diagram, like the Wirtinger presentation.

A ribbon knot is a knot bounding a singular disk with only ribbon singularities. We give here a formal definition of ribbon knot.
\begin{defn}\label{def:ribbon}
Let $r: D^2 \longrightarrow M^3$ be a smooth immersion of an oriented disk into a smooth 3-manifold. If each component of the singularities of $r(D^2)$ is an arc and one component of its preimage in $D^2$ is interior, then $r (\partial D^2 )$ is a \emph{ribbon knot} in $ M^3 $, and $ R = r(D^2)$ is a \emph{ribbon}.
\end{defn}

There were some fine researches on Alexander polynomials for some special types of ribbon knots, especially symmetric unions and simple-ribbon knots. Symmetric unions were first introduced by Kinoshita and Terasaka\cite{KinoshitaTerasaka1957} in 1957 and generalized by Lamm\cite{Lamm2000} in 2000. Lamm obtained many properties of Alexander polynomials of symmetric unions \cite{Lamm2000}. But no one has given a general formula for Alexander polynomials of symmetric unions. Recently, some researchers (c.f.\cite{KishimotoShibuyaTsukamotoIshikawa2021}) defined a special class of ribbon knots in technical terms, called simple-ribbon knots, gave a general formula for their Alexander polynomials, and applied their formula to consider some specific problems.

Ribbon knots are slice knots and Fox conjectured that all slice knots are ribbon knots. 
A classical result proved by Fox and Milnor\cite{FoxMilnor} is that the Alexander polynomial of a slice knot has the form $f(t)f(t^{-1})$ for some integer-coefficient polynomial $f(t)$. 
For ribbon knots, we know that we can give a simpler proof of this result in $S^3$, see \cite{Kauffman1987}. 
A folklore result is that $f(t)$ can be taken to be the Alexander polynomial of the ribbon in $ B^4 $.

In this paper, we consider general ribbon knots, give general formulae for their Alexander polynomials, and study their properties and variants in detail. Furthermore, we always study Alexander polynomial of ribbon knots in a more general setting, i.e., in integer homology spheres.

\subsection{Alexander polynomials of ribbon knots in $ \mathbb{Z}HS^3 $}\label{subsect:introduction1}

We first give a natural definition to record the intrinsic singularity information of an arbitrary ribbon $ R $, called \emph{ribbon diagram}. See Definition \ref{def:rd}. 
We further define \emph{ribbon graph} to record the main information of the ribbon diagram. See Definition \ref{def:rg}. We stress here that ribbon graph is not a graph but a pair,  consisting of a directed tree and a map from its edges to its vertices. 
We then define a definite integer-coefficient polynomial called \emph{half Alexander polynomial}, denoted $ A_R (t) $, which is determined by the ribbon graph. See Definition \ref{def:hap}. 
Our basic result is Theorem \ref{thm:basic} and the explicit formula therein.
\begin{thm}\label{thm:rawalgo0} (Theorem \ref{thm:basic})
Given a ribbon $ R $ for a ribbon knot $ K $ in $ \mathbb{Z}HS^3 $, let $ A_R (t) $ be the half Alexander polynomial read off from the ribbon graph. Then the Conway-normalized Alexander polynomial of $ K $ is $ \Delta (t) = A_R (t)   A_R (t^{-1}) $. 
Especially, Alexander polynomial of a ribbon knot is determined by the ribbon graph of the ribbon.
\end{thm}

We emphasize that, unless the ribbon diagram is rather trivial, each ribbon diagram corresponds to infinitely many different ribbons and thus to infinitely many different ribbon knots. 

Combining more analysis on half Alexander polynomial and ribbon graph, we give a reduced general procedure for computing the Alexander polynomial of ribbon knots, see Subsection \ref{susect:generalprocedure}. Also, we get a nontrivial upper bound of the breath of Alexander polynomial of a ribbon knot in Corollary \ref{cor:breath}. As a direct corollary of Theorem \ref{thm:rawalgo0}, we obtain a result on the determinant of a ribbon knot, i.e., Corollary \ref{cor:determinant}, which generalizes Theorem 2 in \cite{Lamm2000} of Lamm on the determinant of a symmetric union.

For the ribbon knots in $ S^3 $, the last statement of Theorem \ref{thm:rawalgo0}, if we here weaken ribbon graph to ribbon diagram, had been derived implicitly in some literature, e.g. \cite{Terasaka1959}. However, the methods in \cite{Terasaka1959} essentially calculated the fundamental group by diagrams, having nothing to do with Seifert surfaces, so are not valid for ribbon knots in more general settings, such as in $ \mathbb{Z}HS^3 $. 
Furthermore, the formula for Alexander polynomial we give in this theorem, to the best of our knowledge, has not appeared in literature. Even used in the specific case simple-ribbon knots, defined in \cite{KishimotoShibuyaTsukamotoIshikawa2021}, it is not the formula therein.

The middle part of this paper, Section \ref{sec:hapor}, is devoted to exploring in depth the half Alexander polynomials in several ways. 

The main result of this part is two simplified formulae for computing half Alexander polynomials, which are applicable generally. They are the \emph{contracted formula} given in Theorem \ref{thm:algo2} and the \emph{path-type formula} given in Theorem \ref{thm:algo3}. The equivalence of each formula and the original formula in Theorem \ref{thm:basic} is not obvious. Each formula can greatly simplify the computational effort in some common cases. The path-type formula, although applicable only when the ribbon tree is a path, is also general, since each ribbon can be deformed into this shape while keeping the knot type unchanged, see Proposition \ref{prop:path}.

In Section \ref{sec:hapor}, we give a three dimensional topological explanation of half Alexander polynomial in Proposition \ref{prop:nfolklore}, different from the folklore one previously mentioned. 
We also point out with examples that the half Alexander polynomial gives substantially more information about a ribbon than the Alexander polynomial, see Subsection \ref{subsect:example}. 
More importantly, we point out by counterexamples that the half-Alexander polynomial is not an invariant of ribbon knot, even up to multiplication by $\pm t^{\pm n}$ and changing $ t $ by $ t^{-1} $. This has the following implication. In the classical result of Fox and Milnor\cite{FoxMilnor}, although they did not say that $f(t)$ is an invariant of the slice knot, they did not deny the possibility either. While it is certainly possible to have different $f(t)$ when the factorization of the Alexander polynomial is not unique, is it possible that there exists a special factorization such that $f(t)$ is indeed an invariant of slice knot? To the best of our knowledge, no literature has explicitly rejected this possibility. And our counterexample illustrates that this is essentially impossible.

Using our path-type formula, Theorem \ref{thm:algo3}, we revisit a result of Terasaka\cite{Terasaka1959}, that given any integer-coefficient polynomial $ f(t) $ satisfying $ f(1) = \pm 1 $, there exists a ribbon knot whose Alexander polynomial is $ f(t)f(t^{-1}) $. Note that since the factorization of an Alexander polynomial polynomial may not be unique, logically it does not imply the following result.
\begin{cor}\label{thm:converse0} (Corollary \ref{cor:geography})
Given a polynomial $ f(t) $ with integer coefficients and $ f(1) = \pm 1 $, there is a ribbon with fusion number 1 so that the half Alexander polynomial $ A_R (t) \dot{=} f(t) $ .
\end{cor}
However, if we examine Terasaka's construction\cite{Terasaka1959}, we find that his construction already satisfies the above result. We prove the above result using our path-type formula, and our construction is much more compact and regular.

Being a digression, we give a generalization of Theorem \ref{thm:algo2} and the contracted formula, namely Theorem \ref{thm:generalize}, which provides a formula for the Alexander polynomial for the band sum of knots. Although this seems to be of little theoretical significance, it has received a lot of attention. In fact, the whole paper \cite{Terasaka1959} is about computing the Alexander polynomial for the band sum of a knot and an unknot using the Wirtinger representation, and the main point of paper \cite{KishimotoShibuyaTsukamotoIshikawa2021} is to compute the Alexander polynomial for a very special band sum of a knot and many unknots. The method in \cite{KishimotoShibuyaTsukamotoIshikawa2021} also uses Seifert surfaces, but it uses a variety of special techniques such as sliding, winding and tubing, while ours is relatively uniform and general.

\subsection{Alexander polynomials of general knots and virtual knots}\label{subsect:introduction2}

Applying our formula in Section \ref{subsect:halfalexpoly} to a special type of ribbon knot, the connected sum of a knot and its mirror image, unexpectedly inspires us to discover a new formula for computing the Alexander polynomial for a general knot. The application of our formulas to virtual knots sometimes simplifies the calculation even more.

Gauss diagram is a way to represent a knot. Actually, the Gauss diagram represents more general object called virtual knots, and not all the Gauss diagrams represent a knot. Thus it is an interesting problem to explore a Gauss diagram formula for knot invariants, because such a formula may bring an extension of a knot invariant to a virtual knot invariant.

Although Alexander polynomial is a quite classical invariant, looking for new algorithmatic formula for Alexander polynomial is still attracting for some knot theorists in these decades. In the spirit of \cite{GoussarovPolyakViro2000} that any finite type invariant of knots can be computed from a Gauss diagram, by counting suitable sub-Gauss diagrams, in \cite{ChmutovKhouryRossi2009}, the authors gave an explicit Gauss diagram formula for the coefficients of Conway polynomial. Their proof was purely combinatorial, by induction on the number of crossings.

We give a formula for the whole Alexander polynomial in terms of Gauss diagram, which is quite different from the formula of Chmutov, Khoury and Rossi \cite{ChmutovKhouryRossi2009}. Our method is rather topological, by considering infinite cyclic cover space using a proper surface, which is not a Seifert surface. Our formula is generalized to virtual knots and long virtual knots, giving a topological interpretation of Alexander polynomial similar to the classical case.

Although it is reasonable to assume that all diagrammatic formulas for Alexander polynomial should be transformed to each other in certain standard ways, such transformation formulas may be hard to find. We discuss this in Subsection 5.2.2. We explain the difficulty to transform between our formula and the Seifert matrix from canonical Seifert surface and that the reason for this difficulty is actually a benefit, allowing us to generalize our results to virtual knots. We manage to find a transformation formula between our formula and the well-known formula from Wirtinger presentation. Besides, we can derive Chmutov, Khoury and Rossi's formula\cite{ChmutovKhouryRossi2009} for the second coefficient of Conway polynomial from our formula.

We remind again that Subsection 1.1 and Subsection 1.2 are logically independent. However, the parallelism between many of their results indicates an idea that, algebraically, a ribbon knot is like a connected sum of a virtual knot and its mirror image.

\section{Backgrounds}
\subsection{Alexander polynomials of knots in $ \mathbb{Z}HS^3 $}
Let $ \mathbb{Z}HS^3 $ be any integer homology sphere. Linking number is well-defined in $ \mathbb{Z}HS^3 $ exactly as in $S^3$. Let $K$ be a smooth knot in $ \mathbb{Z}HS^3 $. Then $K$ can bound a Seifert surface $F$, i.e., oriented spanning surface. For a basis $[f_1], \dots, [f_{2g}]$ of $H_1(F)$, the Seifert matrix is defined as $A = \left( \text{lk}(f_i, f_j^+) \right)_{2g \times 2g}$, where $f_j^\pm$ denotes the push-off of cycle $f_j$ along the $\pm$ direction of $F$. Then the Conway-normalized Alexander polynomial \cite{Lickorish1997} is
\begin{align}
    \Delta(t) = \left| t^{1/2} A - t^{-1/2} A^T \right|.
\end{align}

\subsection{Alexander polynomials of virtual knots}\label{subsect:virtual Wirtinger}
\begin{figure}[h]
\centering
\includegraphics[scale=0.3]{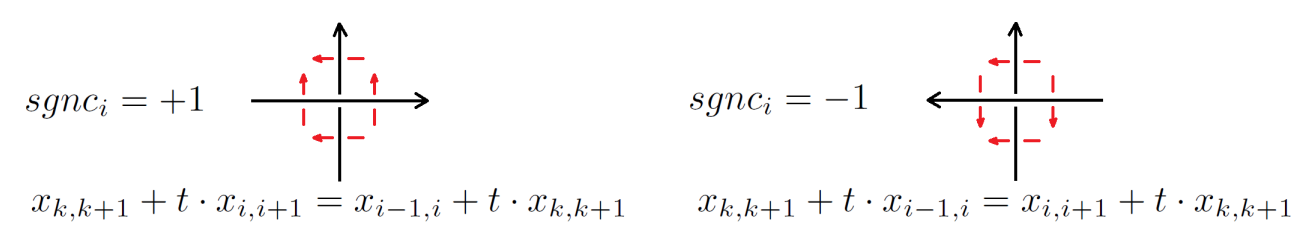}\label{fig:FoxofWirtinger}
\caption{Fox differential of Wirtinger presentation.}
\end{figure}

Let $K$ be a diagram for a virtual knot $K$, with $n$ true crossings $c_1, \dots, c_n$. Orient $K$. For ease of presentation, we label the crossings in order of the under-crossings along $K$. For the over-passing arc between $c_i$ and $c_{i+1}$, denote the standard meridian generator encircling it by $x_{i,i+1}$ (Here $n + 1 = 1$). Let $x_{k(i),k(i)+1}$ be the over-passing arc at each $c_i$. Fig. \ref{fig:FoxofWirtinger} gives an $n \times n$ presentation matrix $W = \Phi - t \Psi$, where $W = \Phi$ and $\Psi$ are two integer matrices. For classical knots, this is the Fox free differential of Wirtinger presentation.

We define the greatest common divisor of the $(n - 1) \times (n - 1)$ minors of $W$ to be the \emph{Alexander polynomial} for the virtual knot, denoted $\Delta_K(t)$.

\begin{prop}
$\Delta_K(t)$ is a virtual knot invariant.
\end{prop}

\begin{proof}
Track the change of $W$ under Reidemeister moves.
\end{proof}

For classical knots, in fact, the determinant of any $(n - 1) \times (n - 1)$ minor of $W$ is the Alexander polynomial. For virtual knots, viewed as knots in thickened surfaces, we can give a topological explanation of $\Delta_K(t)$ as the generator of the homology module of a infinite cyclic cover space, which is constructed similarly as in the proof of Theorem \ref{thm:classical} in Subsection \ref{subsect: classical knots}.

There are many different definitions of Alexander polynomials for virtual knots, defined from different approaches (c.f. \cite{JuhaszKauffmanOgasa2022}). Most of them are 2-variable. To our knowledge, the definition we give here, although a simplest, has not appeared (at least explicitly) in other literature. Since we do not want to go too far in virtual knots here, the relationship between these Alexander polynomials may be discussed in detail in future work.

\section{Alexander module of ribbon knots in $ \mathbb{Z}HS^3 $}

In this section we define ribbon diagram, ribbon graph and half Alexander polynomial for a ribbon and prove Theorem \ref{thm:basic}. After this, we discuss more Alexander invariants of ribbon knot which cannot be determined by the ribbon.

\subsection{Ribbon graph, half Alexander polynomial and our basic theorem}\label{subsect:halfalexpoly}

Recall Definition \ref{def:ribbon}. We shall naturally define a diagram on $D^2$ to represent the singularity information of the ribbon.

\begin{figure}[h]
\centering
\includegraphics[scale=0.33]{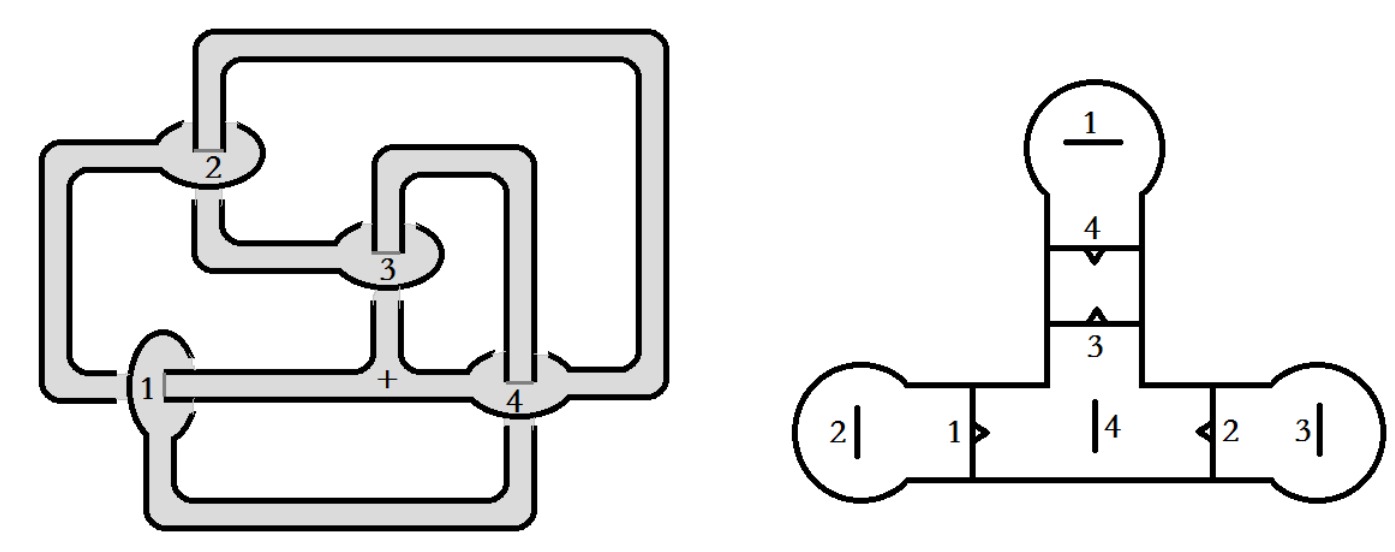}
\caption{A ribbon and its ribbon diagram.}
\label{Fig:a ribbon and its ribbon diagram}
\end{figure}

\begin{defn}\label{def:rd}
    Let $R$ be a ribbon defined in Definition \ref{def:ribbon}. Label the components of singularities of the ribbon by $\alpha_1, \alpha_2, \dots, \alpha_g$. For each $i = 1, \dots, g$, let $\beta_i$ be the interior component of $r^{-1}(\alpha_i)$ in $D^2$ and $\gamma_i$ be the other component of $r^{-1}(\alpha_i)$. If one side of $\gamma_i$ is at the positive normal direction of an embedded neighborhood of $\beta_i$, mark it by a small triangle, and denote the marked $\gamma_i$ by $\hat{\gamma}_i$. The disk with these arcs, $(D^2, \cup_{i=1}^g \beta_i, \cup_{i=1}^g \hat{\gamma}_i)$, is the \textit{ribbon diagram} of $R$.
\end{defn}
See Fig. \ref{Fig:a ribbon and its ribbon diagram} for an example.

We introduce ribbon graph to carry informations of the ribbon diagram.
\begin{defn}\label{def:rg}
Let $(D^2, \cup_{i=1}^g \beta_i, \cup_{i=1}^g \hat{\gamma}_i)$ be a ribbon diagram. The pair $(T, S)$ defined as follows is the \textit{ribbon graph} for the ribbon diagram.

\begin{itemize}
\item \textit{Ribbon tree} $T$: each vertex $v_i$ is a component of $D^2 - \cup_{i=1}^g \gamma_i$ for $i = 1, \dots, g + 1$, and each directed edge $E_i$ corresponds to $\hat{\gamma}_i$ for $i = 1, \dots, g$ so that the head of $E_i$ contains the mark triangle of $\hat{\gamma}_i$ and the tail of $E_i$ is the other component whose boundary contains $\gamma_i$.

\item \textit{Singularity map} $S$: a map from the edge set of $T$ to the vertex set of $T$ such that $S(E_i) = v_j$ if and only if $\beta_i$ is in $v_j$.
\end{itemize}
\end{defn}

Note that each $\hat{\gamma}_i$ separates $D^2$. We point out that there may be different ribbon diagrams corresponding to one ribbon graph. But if $T$ is embedded in a plane, then ribbon graphs correspond bijectively to ribbon diagrams up to homeomorphism.

Now we define a \textit{ribbon matrix} $\rho_{g \times g}$ in terms of the ribbon graph. For each $i = 1, \dots, g$, subdivide $E_i$. Denote the new vertex by $v_{E_i}$ and the obtained directed tree by $T_{E_i}$. Let $P_i$ be the unique path in $T_{E_i}$ from $v_{E_i}$ to $S(E_i)$. Let
\begin{align}
    \rho_{ii} = 
\begin{cases}
\frac{1}{2}, & \text{if the first edge of } P_i \text{ is a forward arc;} \\
-\frac{1}{2}, & \text{if the first edge of } P_i \text{ is a reverse arc.}
\end{cases}
\end{align}
\begin{align}
\rho_{ji} = 
\begin{cases}
1, & \text{if } E_j \text{ is a forward arc of } P_i; \\
-1, & \text{if } E_j \text{ is a reverse arc of } P_i; \\
0, & \text{if } E_j \notin P_i.
\end{cases}
\quad \forall j \neq i. \nonumber
\end{align}\label{for:ribbon matrix rho}

Set matrix $R(t)_{g \times g} = (t - 1)\rho - \frac{1}{2}(t + 1)I$. That is,
\begin{align}
    R_{ii}(t) = 
\begin{cases}
-1, & \text{if the first edge of } P_i \text{ is a forward arc;} \\
-t, & \text{if the first edge of } P_i \text{ is a reverse arc.}
\end{cases}
\end{align}
\begin{align}
    R_{ji}(t) = 
\begin{cases}
t - 1, & \text{if } E_j \text{ is a forward arc of } P_i; \\
1 - t, & \text{if } E_j \text{ is a reverse arc of } P_i; \\
0, & \text{if } E_j \notin P_i.
\end{cases}
\quad \forall j \neq i. \nonumber
\end{align}\label{for:matrix R(t)}

\begin{defn}\label{def:hap}
    The \textit{half Alexander polynomial} of the oriented ribbon $R$ is $A_R(t) = |R(t)|$.
\end{defn}

\begin{lem}
Half Alexander polynomial $A_R(t)$ is an invariant of the oriented ribbon $R$, up to ambient isotopy of $\mathbb{Z}HS^3$.
\end{lem}

\begin{proof}
$R(t)_{g \times g}$ is determined by the ribbon graph, and $|R(t)|$ is independent of the labels of edges in the ribbon tree.
\end{proof}

We remark that the half Alexander polynomial of the ribbon with the other orientation is $(-t)^g |R(t^{-1})|$.

\begin{thm}\label{thm:basic}
Given a ribbon $R$ for a ribbon knot $K$ in $\mathbb{Z}HS^3$, let $A_R(t)$ be the half Alexander polynomial read off from the ribbon graph. 
Then the Conway-normalized Alexander polynomial of $K$ is $\Delta(t) = A_R(t)A_R(t^{-1})$.
\end{thm}

\begin{example}

\begin{figure}[htbp]
    \centering
    \includegraphics[scale=0.3]{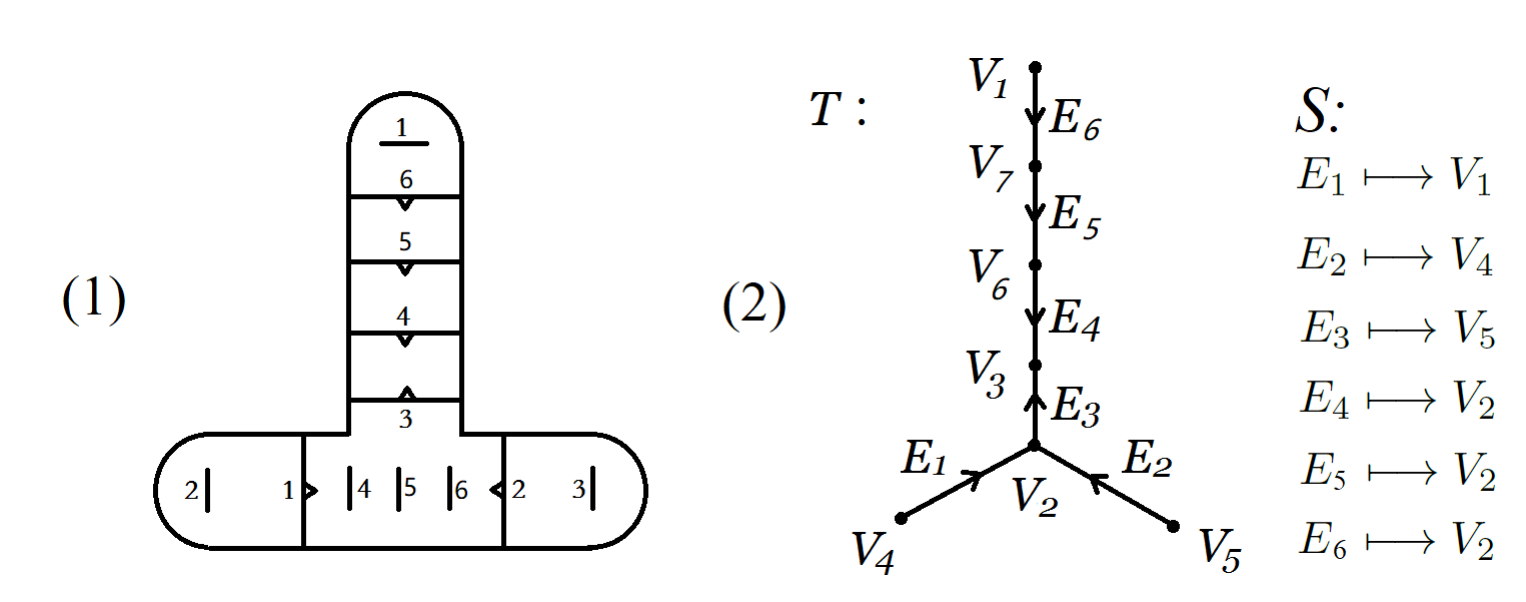}
    \caption{A ribbon diagram and its ribbon graph.}
    \label{fig:ribbon diagram ribbon graph}
\end{figure}

\begin{figure}[htbp]
    \centering
    \includegraphics[scale=0.3]{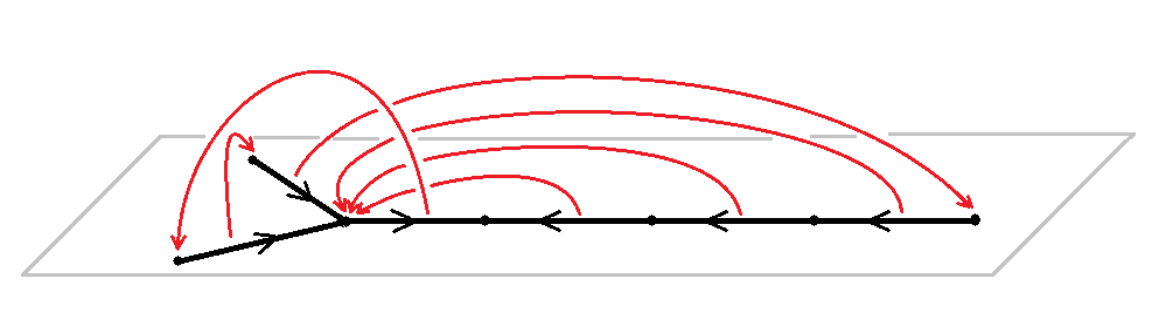}
    \caption{Ribbon graph drawn in a 3 dimensional manner.}
    \label{fig:ribbon graph 3 dimensional}
\end{figure}

For the ribbon in Fig. \ref{fig:ribbon diagram ribbon graph}(1), the ribbon graph $(T, S)$ is as shown in Fig. \ref{fig:ribbon diagram ribbon graph} and also Fig. \ref{fig:ribbon graph 3 dimensional}. By the formula (\ref{for:ribbon matrix rho}), the ribbon matrix should be
\[
\rho = 
\begin{pmatrix}
\frac{1}{2} & -1 & 0 & 0 & 0 & 0 \\
0 & \frac{1}{2} & -1 & 0 & 0 & 0 \\
1 & 0 & -\frac{1}{2} & -1 & -1 & -1 \\
-1 & 0 & 0 & \frac{1}{2} & 1 & 1 \\
-1 & 0 & 0 & 0 & \frac{1}{2} & 1 \\
-1 & 0 & 0 & 0 & 0 & \frac{1}{2}
\end{pmatrix}.
\]
Thus the half Alexander polynomial is
\[
A_R(t) = |(t - 1)\rho - \frac{1}{2}(t + 1)I| = t(1 + t^2 - 3t^3 + 3t^4 - t^5),
\]
and we get the Conway-normalized Alexander polynomial
\[
\Delta(t) = -t^{-5} + 3t^{-4} - 4t^{-3} + 7t^{-2} - 15t^{-1} + 21 - 15t + 7t^2 - 4t^3 + 3t^4 - t^5.
\]
\end{example}

The following corollary generalizes Lamm's Theorem 2.6 in [11] on determinants of symmetric unions.

\begin{cor}\label{cor:determinant}
The determinant of a ribbon knot is independent of the directions of edges in the ribbon tree for any ribbon.
\end{cor}

\begin{proof}
For any ribbon knot $K$ and its ribbon $R$, apply formula (\ref{for:matrix R(t)}) to $\det K = \Delta(-1) = |R(-1)|^2$.
\end{proof}

\subsection{Proof of Theorem \ref{thm:basic}.}

We divide our proof of Theorem \ref{thm:basic} into four steps.

\subsubsection{From ribbon to Seifert surface}\label{sss:From ribbon to Seifert surface}

Let $K$ be a ribbon knot with ribbon $r(D^2)$. We first desingularize the ribbon to get a Seifert surface for $K$. Let $\alpha_i$ be a singularity component, where $i = 1, \dots, g$.

\begin{enumerate}
\item Cut along $\alpha_i$ in the orientation-preserving way and deform the surface locally, which produces a hole and two ribbon-ends.
\item Slide the edge of a ribbon-end along the boundary of the hole towards an arc near the other ribbon-end.
\end{enumerate}

\begin{figure}[htbp]
    \centering
    \includegraphics[scale=0.2]{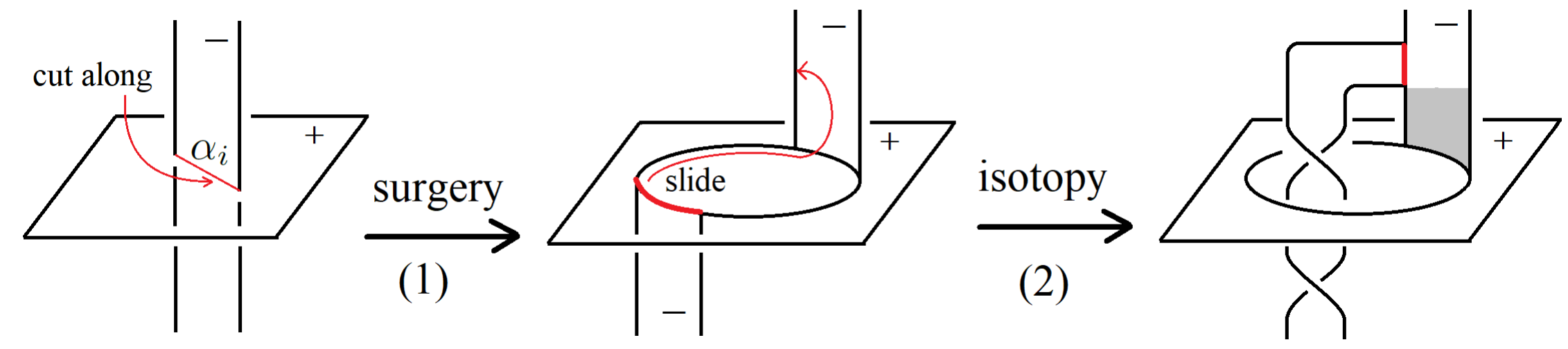}
    \caption{ Surgery and isotopy. }
    \label{fig:surgery and isotopy}
\end{figure}

The resulting surface is a Seifert surface for $K$ of genus $g$, denoted $F_g$. See Fig. \ref{fig:surgery and isotopy} for demonstration.

It is desirable to view $F_g$ from the ribbon diagram. On the ribbon diagram, for each $i = 1, \dots, g$,
\begin{enumerate}
\item Cut along $\beta_i$ to produce a hole, say $\mathring{N}(\beta_i)$.
\item Choose an arc in $\partial N(\beta_i)$ and an arc in $\partial D^2$ containing exactly one endpoint of $\gamma_i$, and glue the two arcs in the orientation-preserving way.
\end{enumerate}

It is easy to see that the surface obtained is $F_g$, as shown in Fig. \ref{fig:a ribbon, a representation on ribbon diagram}(2).

\begin{figure}[htbp]
		  \centering
    \includegraphics[scale=0.23]{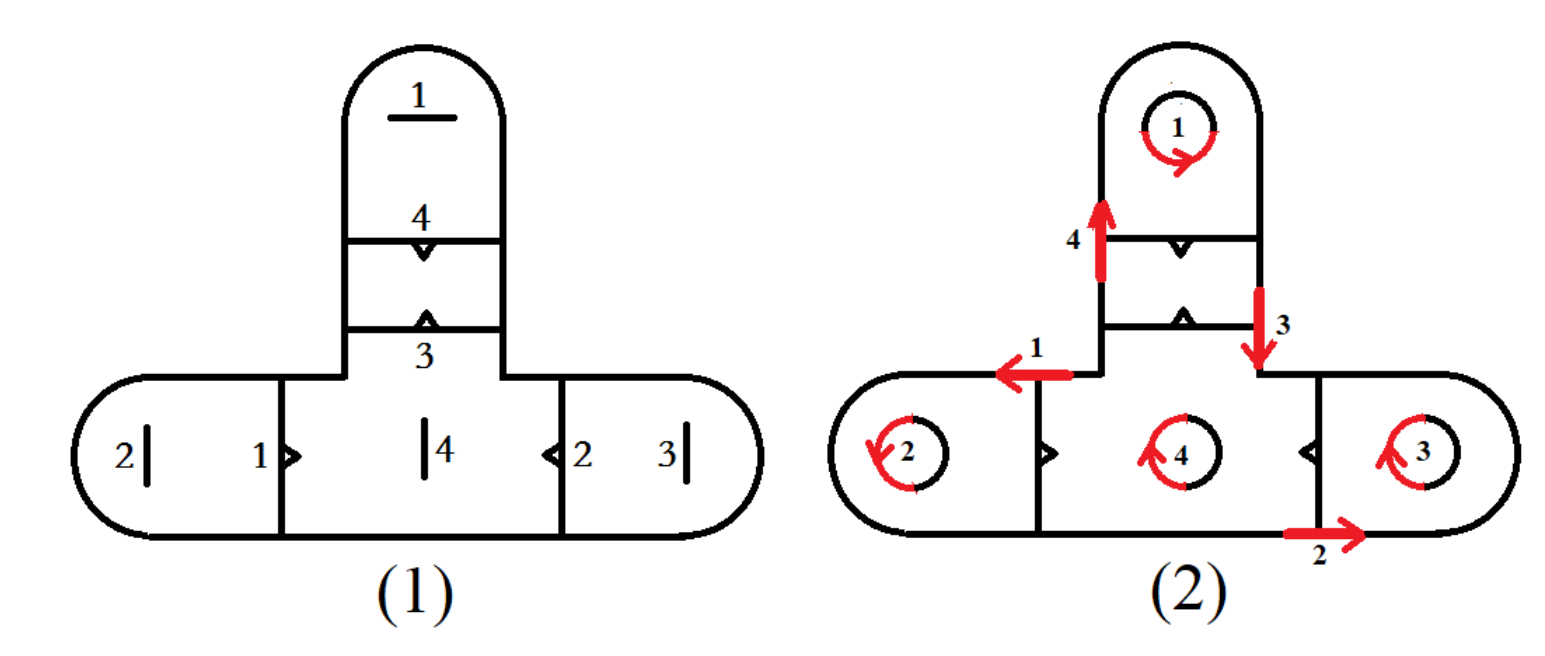}
    \caption{A ribbon, a representation of $F_g$ on ribbon diagram.}
	\label{fig:a ribbon, a representation on ribbon diagram}
\end{figure}

\subsubsection{Homological basis on $F_g$}\label{sss:Homological basis on Fg}

To calculate the Seifert matrix of $F_g$, we first choose a standard basis of $H_1(F_g)$ based on the representation of $F_g$ on the ribbon diagram. For each $i = 1, \dots, g$,
\begin{enumerate}
\item Let $e_i$ be a simple closed curve encircling $\partial N(\beta_i)$, oriented counterclockwise as seen from the positive side of $D^2$.
\item Let $u_i$ be an interior point of $\gamma_i$, and denote the subarc of $\gamma_i$ from the glued endpoint of $\gamma_i$ to $u_i$ by $\gamma_{i0}$.
\item On the ribbon diagram, let $\dot{f}_i$ be a path from $u_i$ to its corresponding point in $\beta_i$. On $F_g$, let $f_i$ be the path composed by $\gamma_{i0}$ and $\dot{f}_i$.
\end{enumerate}

\begin{figure}[htbp]
		  \centering
    \includegraphics[scale=0.23]{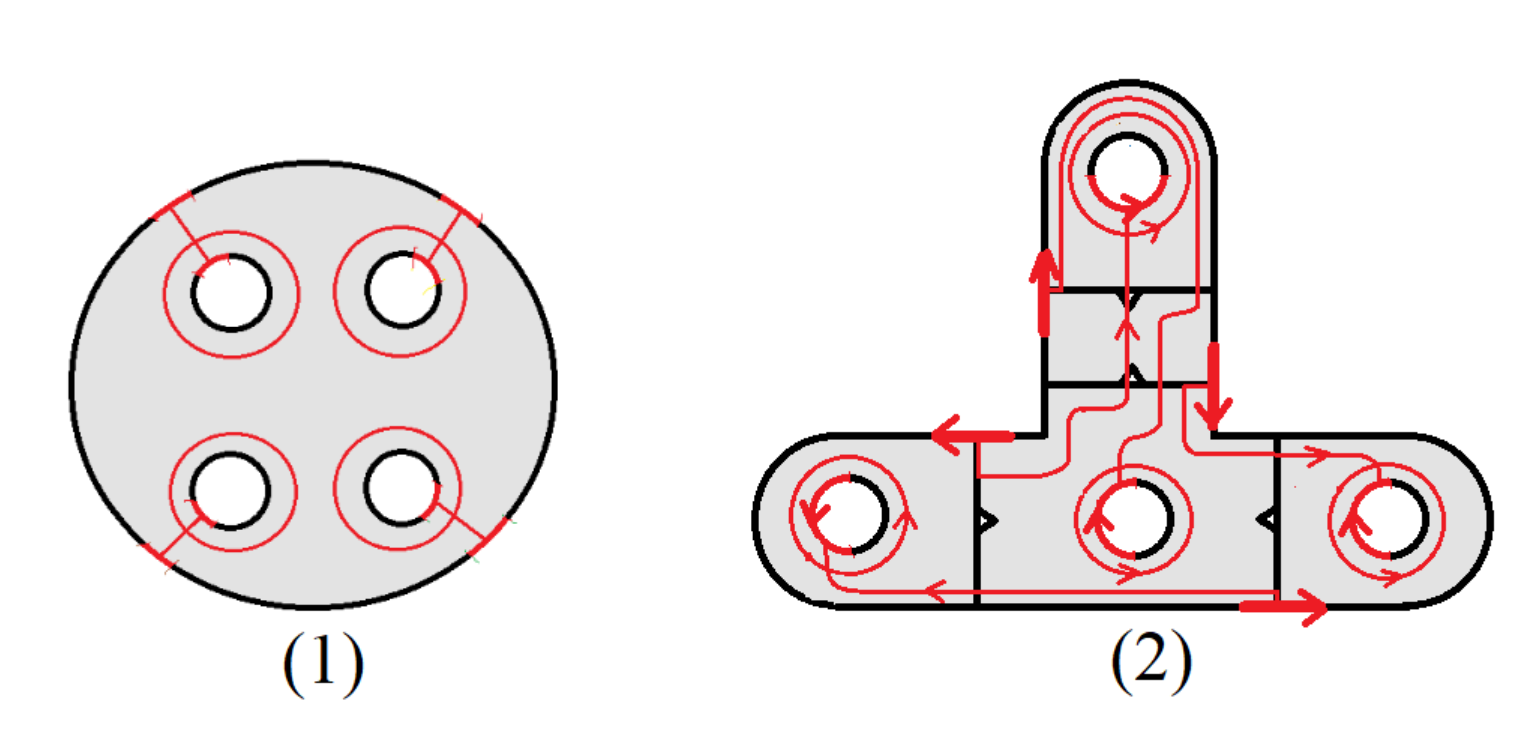}
    \caption{A homological basis on $ F_g$ }
	\label{fig:a homological basis on Fg}
\end{figure}

Choose $f_i$'s mutually disjoint and view them as loops in $F_g$, as demonstrated in Fig. \ref{fig:a homological basis on Fg}(2). Then $[e_1], \dots, [e_g], [f_1], \dots, [f_g]$ form a standard basis of $H_1(F_g)$ with intersection form
\begin{align}
    \begin{pmatrix}
O & I_{g \times g} \\
-I_{g \times g} & O
\end{pmatrix}.
\end{align}

\begin{rem}
A more accurate way to view $F_g$ and the homological basis from the ribbon is as illustrated in Fig. \ref{fig:glue along the blue arc.}. The boundary of the gray disk contacts an endpoint of $\beta_i$. After desingularization, the gray disk is dug and deformed to a sector. We glue it between the severed $\gamma_i$ to get Fig. \ref{fig:a homological basis on Fg}(2).

\begin{figure}[htbp]
		  \centering
    \includegraphics[scale=0.2]{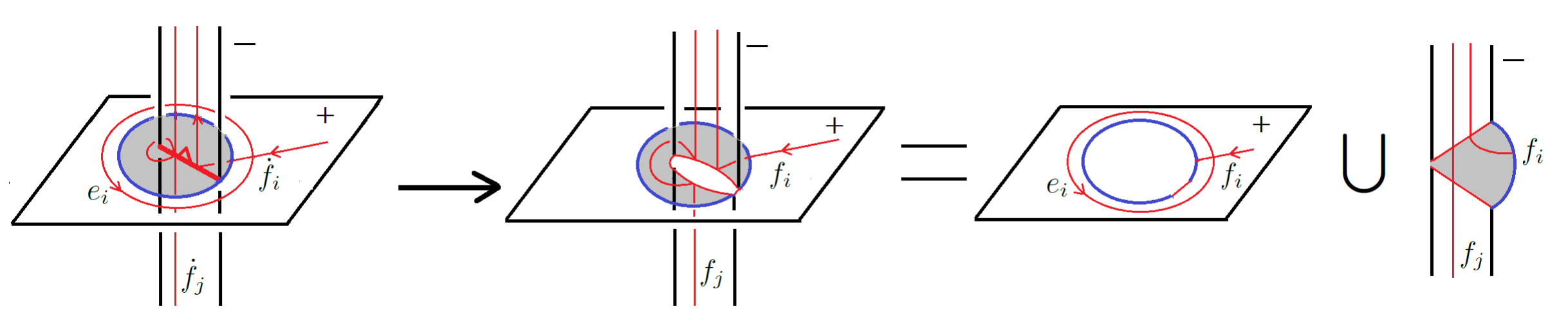}
    \caption{Glue along the blue arc.}
	\label{fig:glue along the blue arc.}
\end{figure}
\end{rem}

\subsubsection{Seifert matrix}\label{sss:Seifert matrix}

Recall the definition of the Seifert matrix (c.f. [10, 11]). Choose the basis $[e_1], \dots, [e_g], [f_1], \dots, [f_g]$ of $H_1(F_g)$ as above. Notice that $\text{lk}(e_i, e_j^+) = 0$ for any $i, j = 1, \dots, g$. The Seifert matrix of $F_g$ has the form
\begin{align}
    A = 
\begin{pmatrix}
O & P \\
Q & L
\end{pmatrix},
\end{align}\label{for:Seifert matrix A}
where
\[
P = (\text{lk}(e_i, f_j^+))_{g \times g}, \quad Q = (\text{lk}(f_i, e_j^+))_{g \times g}, \quad L = (\text{lk}(f_i, f_j^+))_{g \times g}.
\]

Our main task is to calculate $P$ and $Q$.

\begin{figure}[htbp]
		  \centering
    \includegraphics[scale=0.2]{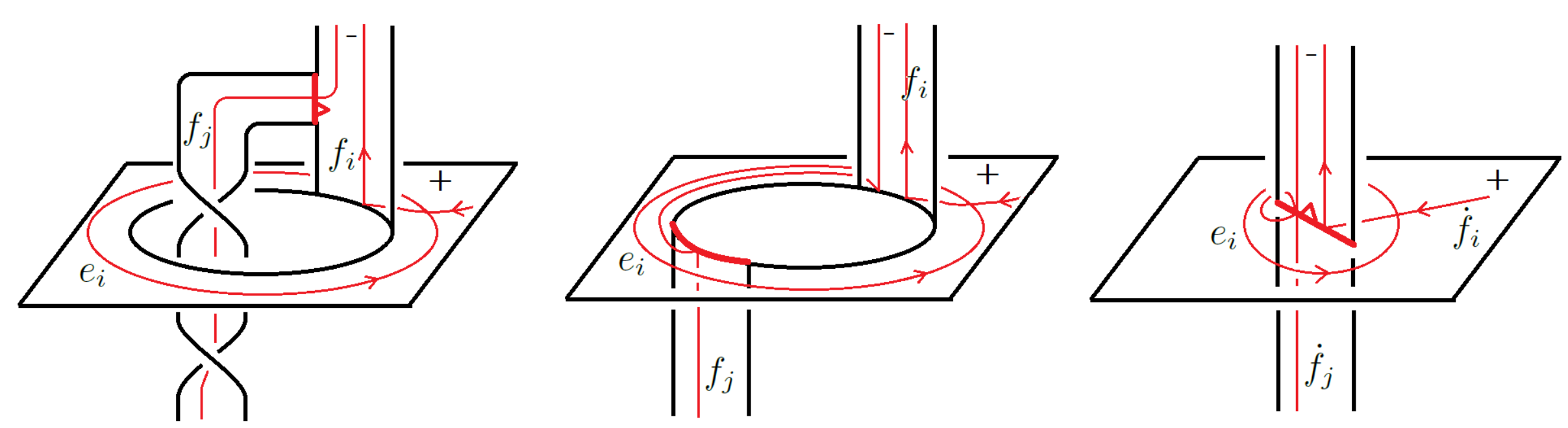}
    \caption{Homological basis locally represented on $F_g$ and on the ribbon.}
	\label{fig:homological basis locally represented}
\end{figure}

Near each singularity $\alpha_i$, we depict $f_i$ and an arbitrarily other $f_j$ on $F_g$ in Fig. \ref{fig:homological basis locally represented}. We have the following observations.

For any $i \neq j$,
\begin{align}
    \text{lk}(e_i, f_j^+) = \text{lk}(e_i^+, f_j) = \text{lk}(e_i, \dot{f}_j),
\end{align}
where $\dot{f}_j$ is viewed as a loop on the ribbon. Consider a point in $\dot{f}_j \cap \gamma_i$. If $\dot{f}_j$ passes through $\gamma_i$ at this point in the same direction as the mark of $\gamma_i$, take the sign at this point as $+1$; if the direction is opposite to the mark of $\gamma_i$, take the sign as $-1$. Then $\text{lk}(e_i, \dot{f}_j)$ equals the sum of signs at the intersection points of $\dot{f}_j \cap \gamma_i$. Notice that each $\gamma_i$ separates $D^2$. As a result, we obtain
\begin{align}
    \text{lk}(e_i, \dot{f}_j) = 
\begin{cases}
1, & \text{if } \beta_j \text{ belongs to the marked side of } \hat{\gamma}_i, \text{ but } \gamma_j \text{ is not;} \\
-1, & \text{if } \gamma_j \text{ belongs to the marked side of } \hat{\gamma}_i, \text{ but } \beta_j \text{ is not;} \\
0, & \text{if } \beta_j \text{ and } \gamma_j \text{ belong to the same side of } \hat{\gamma}_i.
\end{cases}
\end{align}

For $i = 1, \dots, g$, a similar argument, the details of which we omit, suggests that
\[
\text{lk}(e_i, f_i^+) = 0, \quad \text{lk}(f_i, e_i^+) = 1
\]
if $\beta_i$ belongs to the marked side of $\hat{\gamma}_i$;
\[
\text{lk}(e_i, f_i^+) = -1, \quad \text{lk}(f_i, e_i^+) = 0
\]
if $\beta_i$ belongs to the unmarked side of $\hat{\gamma}_i$.

\subsubsection{Read from ribbon graph}

Recall that given a Seifert matrix $A$, the Conway-normalized Alexander polynomial [13] is defined as
\begin{align}
    \Delta(t) = |t^{\frac{1}{2}}A - t^{-\frac{1}{2}}A^T|. \label{for:Alexander definition}
\end{align}

Substituting (\ref{for:Seifert matrix A}) into (\ref{for:Alexander definition}) gives
\begin{align}
    \Delta(t) = |t^{\frac{1}{2}}A - t^{-\frac{1}{2}}A^T| = |tP - Q^T|.
\end{align}

Reading $P$ and $Q$ from the ribbon graph $(T, S)$, we then obtain
\begin{align}
    P = \rho - \frac{1}{2}I, \quad Q = \rho^T + \frac{1}{2}I.
\end{align}
Therefore $tP - Q^T = R(t)$. This completes the proof of Theorem \ref{thm:basic}.

\subsection{Alexander invariants of ribbon knots}

In this subsection, we focus on information other than the Alexander polynomial obtained from the Alexander module. By (\ref{for:Seifert matrix A}), we have a presentation matrix of the Alexander module for the ribbon knot
\begin{align}
   tA - A^T = 
\begin{pmatrix}
O & tP - Q^T \\
tQ - P^T & (t - 1)L
\end{pmatrix} = 
\begin{pmatrix}
O & R(t) \\
-tR(t^{-1})^T & (t - 1)L
\end{pmatrix}. \label{for:tAAT}   
\end{align}

According to \cite{Bai2023}, if we sacrifice the symmetry of matrix $L$, $L$ can be replaced by $\dot{L} \pm QP$, where $\dot{L} = (\text{lk}(\dot{f}_i, \dot{f}_j^+))_{g \times g}$ is symmetric. It is easy to convince that by twining and twisting the partial bands of a ribbon (changing the knot type), one can take any symmetric integer matrix to be $\dot{L}$. Similarly, $L$ in (\ref{for:tAAT}) can be an arbitrary symmetric integer matrix. As the ribbon diagram only determines $R(t)$, ribbon knots with the same ribbon diagram may have different Alexander invariants.

\subsubsection{Ribbon diagram does not determine Alexander invariants.}

It is known that the ribbon knots $6_1$ and $9_{46}$ have the same Alexander polynomial but different second elementary ideals (c.f. [13]).

\begin{figure}[htbp]
		  \centering
    \includegraphics[scale=0.25]{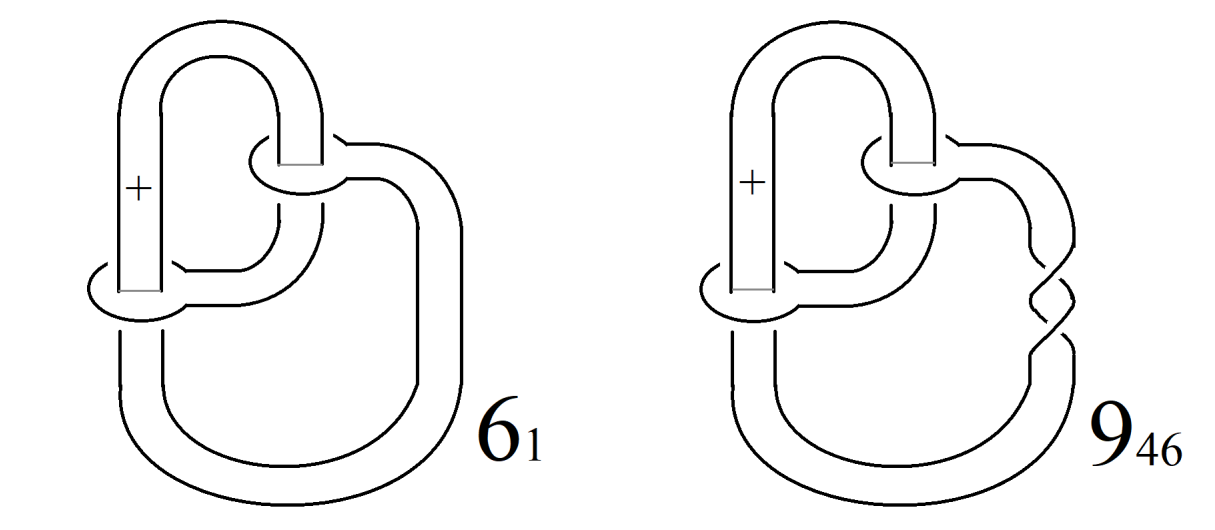}
    \caption{Ribbon for $6_1$ and ribbon for $9_46$.}
	\label{fig:10}
\end{figure}

\begin{example}
If we use the ribbons for $6_1$ and for $9_{46}$ given in Fig. 10, we get two $tA - A^T$ only different in the block $(t - 1)L$,
\[
\begin{pmatrix}
0 & 0 & -1 & 1 - t \\
0 & 0 & 1 - t & -1 \\
t & 1 - t & 1 - t & 0 \\
1 - t & t & 0 & 0
\end{pmatrix}_{n \times n}, \quad 
\begin{pmatrix}
0 & 0 & -1 & 1 - t \\
0 & 0 & 1 - t & -1 \\
t & 1 - t & 0 & 0 \\
1 - t & t & 0 & 0
\end{pmatrix}.
\]
We can calculate that $\Delta_{6_1}(t) = \Delta_{9_{46}}(t) = 2t^{-1} - 5 + 2t$, but the second Alexander ideals of $6_1$ and $9_{46}$ are $\langle 1 \rangle$ and $\langle t + 1, 3 \rangle$ respectively.
\end{example}

\subsubsection{Blanchfield pairing}

Blanchfield pairing is a sesquilinear pairing on the Alexander module which takes values in $\mathbb{Q}(t)/\mathbb{Z}[t, t^{-1}]$. According to [5], under the dual basis of a homological basis of a Seifert surface, a presentation matrix of the Blanchfield pairing is $(t - 1)(A - tA^T)^{-1}$, where $A$ is the Seifert matrix. Therefore, for ribbon knots, using (\ref{for:Alexander definition}), after calculation, we can present the Blanchfield pairing in terms of $R(t)$ and $L$ as
\[
\mathbb{Z}[t, t^{-1}]^{2g}/(tA - A^T) \times \mathbb{Z}[t, t^{-1}]^{2g}/(tA - A^T) \longrightarrow \mathbb{Q}(t)/\mathbb{Z}[t, t^{-1}],
\]
\[
(v, w) \mapsto v^T(t - 1)
\begin{pmatrix}
-R(t)^{-T} & O \\
O & I
\end{pmatrix}
\begin{pmatrix}
(t - 1)L & I \\
I & O
\end{pmatrix}
\begin{pmatrix}
t^{-1}R(t^{-1})^{-1} & O \\
O & I
\end{pmatrix}
\bar{w},
\]
which also depends on $L$.

\section{The half Alexander polynomial of ribbons.}\label{sec:hapor}

Recall that half Alexander polynomial $ A_{R}(t)$ is a definite integer-coefficient polynomial defined for an oriented ribbon $R$, and the Alexander polynomial of the ribbon knot is $ A_R (t)   A_R (t^{-1}) $
We investigate half Alexander polynomials in depth in this section, especially giving two simplified formulae to compute half Alexander polynomial, that is, Theorem \ref{thm:algo2} and Theorem \ref{thm:algo3}.

\subsection{Reductions and procedure for computing Alexander polynomial.}\label{susect:generalprocedure}
We introduce four transformations for ribbon graphs. Let $(T,S)$ be a ribbon graph defined as above.

\begin{figure}
  \begin{minipage}[t]{0.5\linewidth}
    \centering
    \includegraphics[scale=0.24]{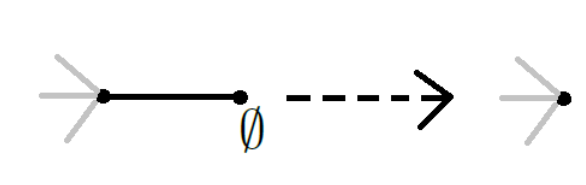}
    \caption{R0-reduction}
    \label{fig:r0-reduction}
  \end{minipage}%
  \begin{minipage}[t]{0.5\linewidth}
    \centering
    \includegraphics[scale=0.24]{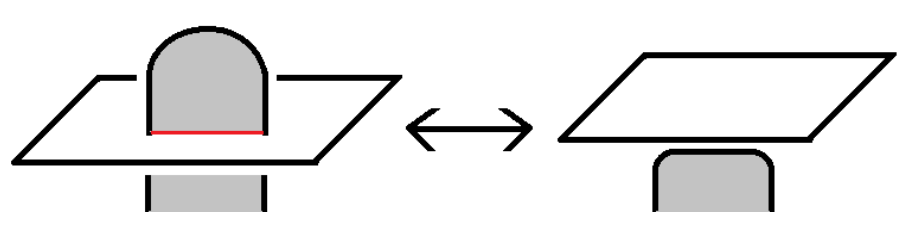}
    \caption{Ribbon 0-move.}
    \label{fig:ribbon 0-move}
  \end{minipage}
\end{figure}

\begin{figure}
  \begin{minipage}[t]{0.5\linewidth}
    \centering
    \includegraphics[scale=0.24]{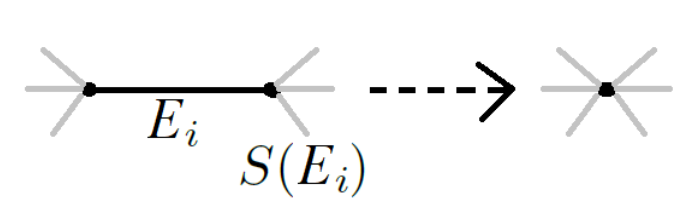}
    \caption{R1-reduction.}
    \label{fig:r1-reduction}
  \end{minipage}%
  \begin{minipage}[t]{0.5\linewidth}
    \centering
    \includegraphics[scale=0.24]{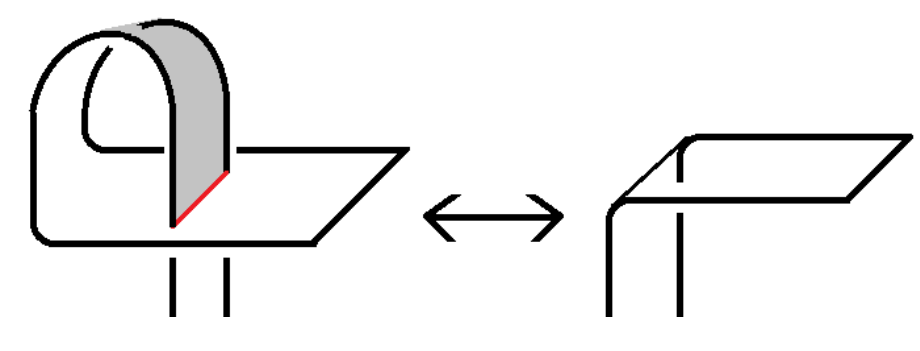}
    \caption{Ribbon 1-move.}
    \label{fig:ribbon 1-move}
  \end{minipage}
\end{figure}

\begin{figure}
  \begin{minipage}[t]{0.5\linewidth}
    \centering
    \includegraphics[scale=0.24]{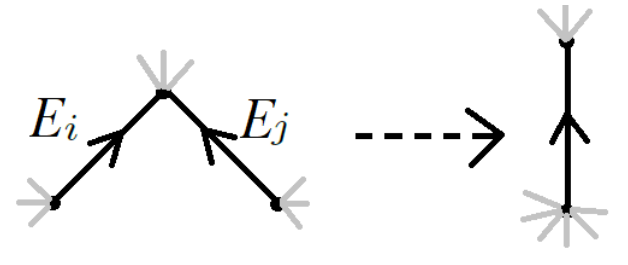}
    \caption{F-reduction.}
    \label{fig:f-reduction}
  \end{minipage}%
  \begin{minipage}[t]{0.5\linewidth}
    \centering
    \includegraphics[scale=0.24]{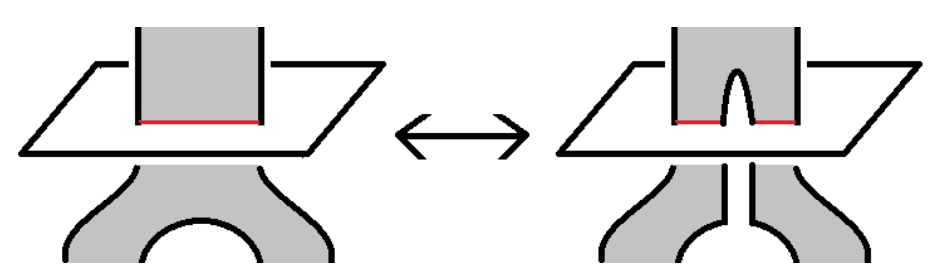}
    \caption{Finger move.}
    \label{fig:finger move}
  \end{minipage}
\end{figure}

\emph{R0-reduction:} If $v_{i}$ is a pendant vertex of $T$ and $v_{i}\notin \text{Im}(S)$, delete $v_{i}$ from $T$ to get a new directed tree, as shown in Fig. \ref{fig:r0-reduction}. Let the new singularity map be the original $S$ restricted to the edges of the new directed tree.

\emph{R1-reduction:} If $S(E_{i})$ is an end of $E_{i}$, contract $E_{i}$ to get a new directed tree, as shown in Fig. \ref{fig:r1-reduction}. Let the new singularity map be $S$ restricted to the edges of the new directed tree and if an edge was mapped to an end of $E_{i}$, we now map it to the new vertex.

\emph{F-reduction:} Suppose $E_{i}$ and $E_{j}$ have the common head or the common tail, and $S(E_{i})=S(E_{j})$. Identify the other ends of $E_{i}$ and $E_{j}$ to get a new directed tree. If an edge was mapped to one of these two ends, we now map it to the new vertex. Fig. \ref{fig:f-reduction} shows the change of $T$ when $E_{i}$ and $E_{j}$ have the common head.

\begin{figure}[htbp]
		  \centering
    \includegraphics[scale=0.25]{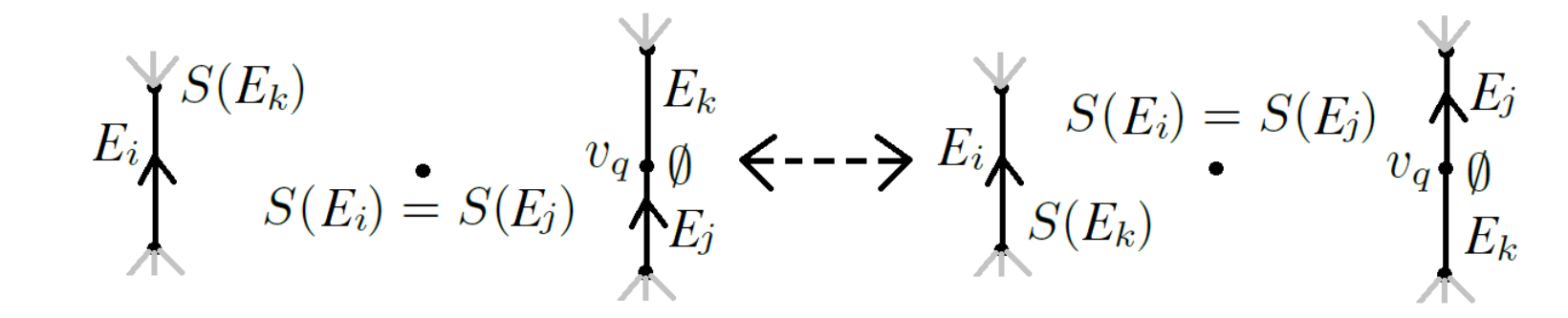}
    \caption{R3-transformation.}
	\label{fig:r3-reduction}
\end{figure}

\emph{R3-transformation:} Suppose $v_{q}$ is of degree 2 with $v_{q}\notin \text{Im}(S)$, $v_{q}$ is the head of $E_{j}$, the other edge incident to $v_{q}$ is $E_{k}$, and $S(E_{i})=S(E_{j})$ where the head of $E_{i}$ is $S(E_{k})$. Interchange $E_{j}$ and $E_{k}$ to get a new directed tree. Revise the singularity map by mapping $E_{k}$ to the tail of $E_{i}$. Such transformation of $(S,T)$ and its reverse are both called R3-transformation. See Fig. \ref{fig:r3-reduction}.

\begin{figure}[htbp]
		  \centering
    \includegraphics[scale=0.23]{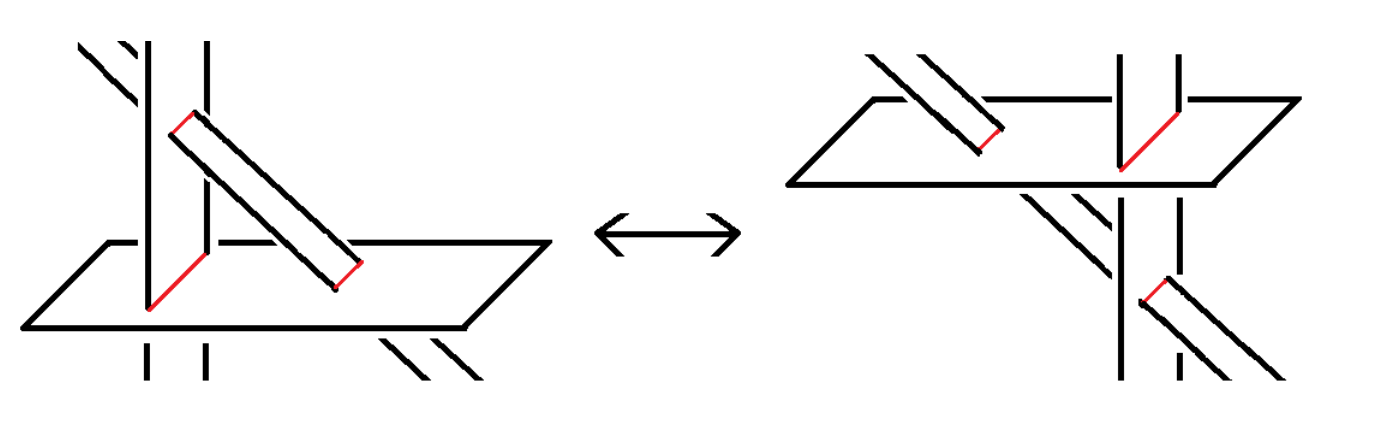}
    \caption{Ribbon 3-move.}
	\label{fig:ribbon 3-move}
\end{figure}

\begin{figure}[htbp]
		  \centering
    \includegraphics[scale=0.23]{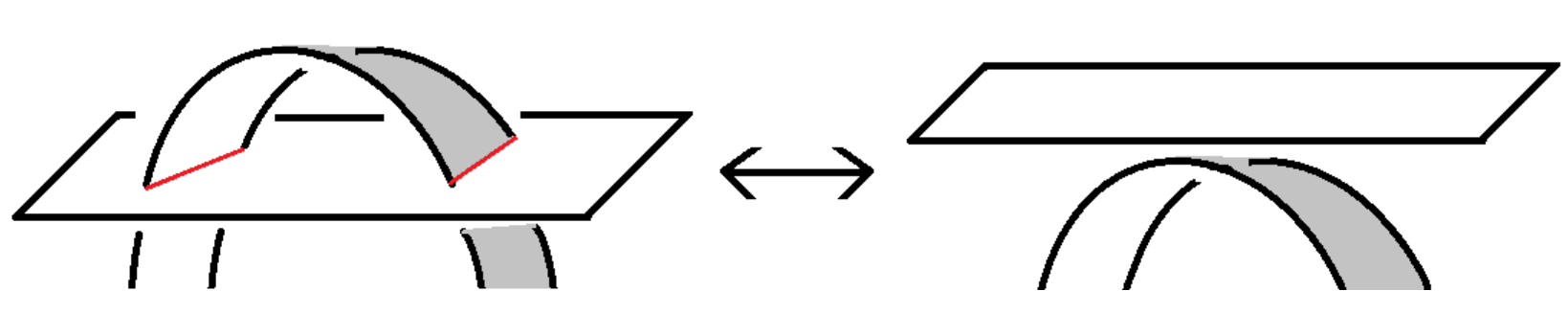}
    \caption{Ribbon 2-move.}
	\label{fig:ribbon 2-move}
\end{figure}

\begin{figure}[htbp]
		  \centering
    \includegraphics[scale=0.21]{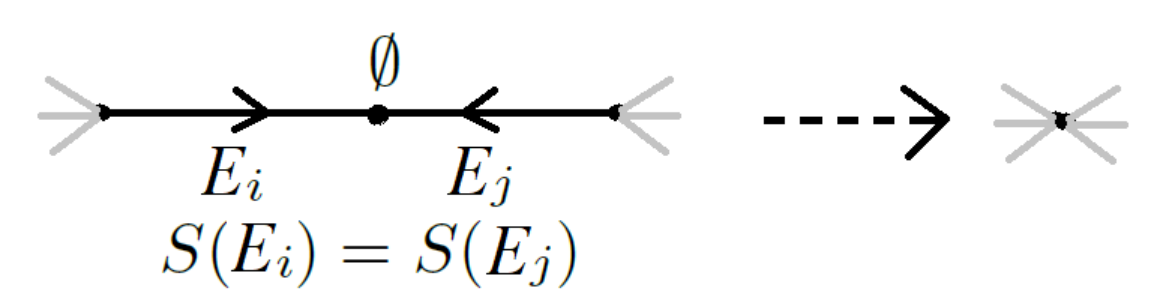}
    \caption{R2-reduction.}
	\label{fig:r2-reduction}
\end{figure}

\begin{lem}\label{lem:reductions unchange}
Under each of the above transformations, the half Alexander polynomial $A_{R}(t)$ changes by a multiplication with $\pm 1$ or $\pm t^{\pm 1}$, and thus the Conway-normalized Alexander polynomial keeps unchanged.
\end{lem}

\begin{proof}
Verify by tracking the change of ribbon matrix.
\end{proof}

We also define several transformations of the ribbon without changing the ribbon knot type. They are \emph{ribbon 0-move} as illustrated in Fig. \ref{fig:ribbon 0-move}, \emph{ribbon 1-move} as illustrated in Fig. \ref{fig:ribbon 1-move}, \emph{finger move} as illustrated in Fig. \ref{fig:finger move} and \emph{ribbon 3-move} as illustrated in Fig. \ref{fig:ribbon 3-move}. It is straightforward to show that these moves induce R0-, R1-, F-reductions and R3-transformation on the ribbon graph respectively.

\begin{rem}
The Ribbon 2-move as shown in Fig. \ref{fig:ribbon 1-move} can be easily viewed as a combination of a ribbon 0-move and a finger move.
Similarly, the R2-reduction as defined in Fig. \ref{fig:r2-reduction} is a combination of R0-reduction and F-reduction.
\end{rem}

R0-reduction of the ribbon graph corresponds exactly to ribbon 0-move of the ribbon. However, we point out that R1-, F-reduction and R3-transformation are not necessarily obtained from ribbon 1-move, finger move, and ribbon 3-move respectively. Therefore the following corollary is not trivial.

\begin{cor}\label{cor:breath}
For a ribbon knot $K$, the breath of $\Delta(t)$ is no bigger than twice the number of edges of $T$ after R0-, R1-, F-reductions and R3-transformations for any ribbon for $K$.
\end{cor}

We now give a general procedure for computing Conway-normalized Alexander polynomial of a ribbon knot from its ribbon diagram. Let $(D^{2},\cup_{i=1}^{g}\beta_{i},\cup_{i=1}^{g}\tilde{\gamma}_{i})$ be a ribbon diagram.

\textbf{Input:} Ribbon graph $(T,S)$.

\begin{enumerate}
    \item Perform R0-, R1-, and F-reductions on $(T,S)$;
    \item Use formula (\ref{for:ribbon matrix rho}) (or (\ref{for:contracted Rijt}) or (\ref{for:path type}) given later) to get the matrix $\rho$;
    \item Compute the polynomial $A_{R}(t)=|(t-1)\rho-\frac{1}{2}(t+1)I|$.
\end{enumerate}

\textbf{Output:} Conway-normalized Alexander polynomial $\Delta(t)=A_{R}(t)A_{R}(t^{-1})$.

\begin{rem}
R3-transformation does not reduce $(T,S)$ directly, but sometimes produces new $(T,S)$ where R0-, R1-, or R2-reductions can be performed.
\end{rem}

\subsection{Contracted formula and Theorem \ref{thm:algo2}.}\label{subsect:Contracted formula}
\begin{figure}
    \centering
    \includegraphics[scale=0.27]{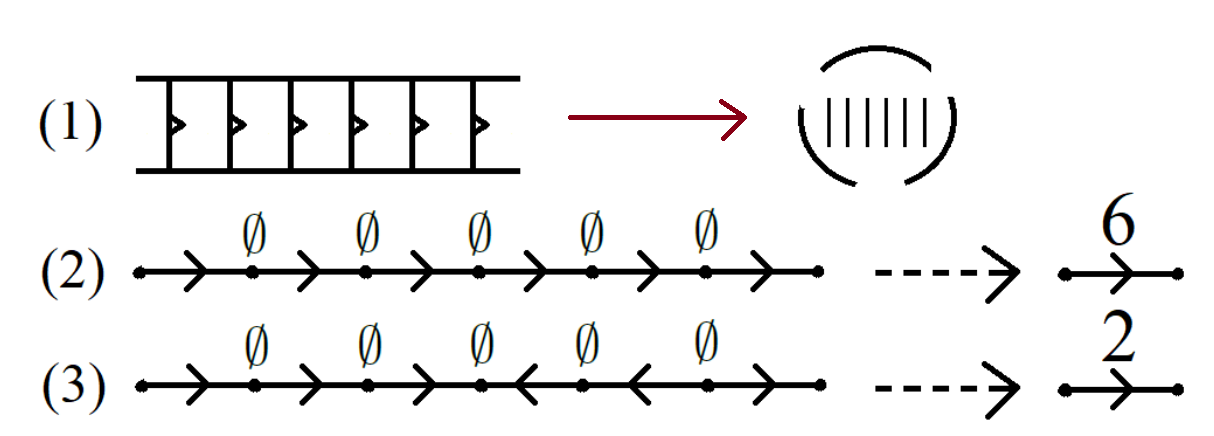}
    \caption{Contracted ribbon graph.}
    \label{fig:contracted ribbon graph}
\end{figure}

The local configurations as in Fig. \ref{fig:contracted ribbon graph}(1) frequently occur in a ribbon diagram, so it is convenient to define a contracted ribbon graph as in Fig. \ref{fig:contracted ribbon graph}(2).

\begin{defn}
Let $(T,S)$ be a ribbon graph. The pair $(T^{\omega},S)$ defined as follows is the \emph{contracted ribbon graph}.

\begin{itemize}
    \item \emph{Weighted tree} $T^{\omega}$: For each maximal induced directed path $P$ in $T$ whose edges are mapped by $S$ to the same vertex and whose interior vertices are not in $\text{Im}(S)$, replace $P$ in $T$ by an edge $E$ with the same direction and endpoints as $P$, and label $E$ with the length of $P$ as its \emph{weight}, denoted $\omega(E)$.
    \item \emph{Singularity map} $S$: Map each edge of $T^{\omega}$ to the image of its original path in $T$.
\end{itemize}
\end{defn}

Using R2-reduction in Fig. \ref{fig:r2-reduction}, to compute Alexander polynomial, we can actually further contract the ribbon graph by replacing``maximal induced directed path'' in the definition of weighted tree by ``maximal induced path'', as shown in Fig. \ref{fig:contracted ribbon graph}(3). 

Using contracted ribbon graph, we give a variant formula to compute half Alexander polynomial of ribbons, usually simpler in practice.

We define a matrix $\boldsymbol{R}(t)$ from a contracted ribbon graph. For each edge $E_{i}$ in $T^{\omega}$, subdivide it. Denote the new vertex by $v_{E_{i}}$ and the obtained directed tree by $T_{E_{i}}$. Let $P_{i}$ be the unique path in $T_{E_{i}}$ from $v_{E_{i}}$ to $S(E_{i})$. Let
\begin{align}
\boldsymbol{R}_{ii}(t) & = 
\begin{cases}
-1, & \text{if the first edge of } P_{i} \text{ is a forward arc}; \\
-t, & \text{if the first edge of } P_{i} \text{ is a reverse arc}. \\
\end{cases} \nonumber \\
\boldsymbol{R}_{ji}(t) & = 
\begin{cases}
t^{\omega_{i}}-1, & \text{if } E_{j} \text{ is a forward arc of } P_{i}; \\
1-t^{\omega_{i}}, & \text{if } E_{j} \text{ is a reverse arc of } P_{i}; \\
0, & \text{if } E_{j}\notin P_{i}.
\end{cases} \quad \forall j\neq i, \label{for:contracted Rijt}
\end{align}
where $\omega_{i}$ is the weight of $E_{i}$.

\begin{thm}\label{thm:algo2}
The half Alexander polynomial $A_{R}(t)=\pm |( {\textbf{R}}_{ji} (t) )| $.
\end{thm}

\begin{proof}
We begin by considering ribbon graph $(T,S)$. After relabelling, for simplification of notations, we may assume $\hat{\gamma}_{1},\hat{\gamma}_{2},...,\hat{\gamma}_{\omega}$ are successive arcs in the ribbon diagram as in Fig. \ref{fig:contracted ribbon graph}(1) so that the path $E_{1}E_{2}...E_{\omega}$ in $T$ corresponds to an edge $E$ in $T^{\omega}$ with $\omega=\omega(E)$. We may assume $\hat{\gamma}_{2},...,\hat{\gamma}_{\omega}$ are in the path $P_{1}$ as defined before (\ref{for:ribbon matrix rho}) in Subsection \ref{subsect:halfalexpoly}.

Suppose $\hat{\gamma}_{1}$ is a forward arc of $P_{1}$. By (\ref{for:matrix R(t)}), the minor

\begin{equation}
(R_{ji}(t))_{1,...,\omega;1,...,\omega}=\left(\begin{array}[]{ccccc}
-1 & & & & \\
t-1 & -1 & & & \\
t-1 & t-1 & -1 & & \\
\vdots & \vdots & \vdots & \ddots & \\
t-1 & t-1 & t-1 & \cdots & -1
\end{array}\right).
\end{equation}

It is easy to see from ribbon graph that

\begin{align}
R_{ji}(t) &=R_{j,\omega}(t), &\quad \forall j\neq 1,2,...,\omega, \quad \forall i=1,2,...,\omega; \nonumber \\
R_{ji}(t) &=R_{j,\omega+1}(t), &\quad \forall j=1,2,...,\omega, \quad \forall i\neq 1,2,...,\omega. \label{for:Rijt 1tow rows}
\end{align}

For $i=1,2,...,\omega-1$, subtract the $(i+1)$th column from the $i$th column, which means choosing $[f_{1}]-[f_{2}],...,[f_{\omega-1}]-[f_{\omega}]$ instead of $[f_{1}],...[f_{\omega-1}]$ in the basis. Then the minor becomes

\begin{equation}
(R_{ji}(t))_{1,...,\omega;1,...,\omega}=\left(\begin{array}[]{cccccc}
-1 & & & & \\
t & -1 & & & \\
0 & t & -1 & & \\
\vdots & \vdots & \vdots & \ddots & \\
0 & 0 & 0 & & -1 \\
0 & 0 & 0 & \cdots & t & -1
\end{array}\right)
\end{equation}
and now

\begin{equation}
R_{ji}(t)=0, \quad \forall j\neq 1,2,...,\omega, \quad \forall i=1,2,...,\omega-1. \nonumber
\end{equation}

Multiply the $1$st, $2$nd ,..., $(\omega-1)$th rows by $t^{\omega-1},t^{\omega-2},...,t$ respectively and add together to the $\omega$th row. Then the new $R_{\omega,\omega-1}(t)$ is $0$. By (\ref{for:Rijt 1tow rows}), for $i\neq 1,2,...,\omega$, the new $R_{\omega,i}(t)$ is $(1+t+\cdots+t^{\omega-1})$ times the old $R_{\omega,i}(t)$. Now we see

\begin{equation}
|\left(R_{ji}(t)\right)|=(1+t+\cdots+t^{\omega-1})\left|\begin{array}[]{cccccc}
-1 & & & & 0 & \\
t & -1 & & & \vdots & *\\
\vdots & \vdots & \ddots & & \vdots & \\
0 & 0 & \cdots & -1 & 0 &  \\
0 & 0 & \cdots & 0  & \frac{-1}{1+t+\cdots+t^{\omega-1}}  & * \\
O &   &        & *  & * & *
\end{array}\right|
\end{equation}
where $\frac{-1}{1+t+\cdots+t^{\omega-1}}$ is the $(\omega,\omega)-$entry. Thus

\begin{equation}
|\left(R_{ji}(t)\right)|=(-1)^{\omega-1}(1+t+\cdots+t^{\omega-1})|\left(R^{*}_{ji}(t)\right)|,
\end{equation}
where $(R^{*}_{ji}(t))$ is the minor of the original $(R_{ji}(t))$ deleting the $1$st, $2$nd ,..., $(\omega-1)$th rows and the $1$st, $2$nd ,..., $(\omega-1)$th columns and replacing $-1$ by $\frac{-1}{1+t+\cdots+t^{\omega-1}}$ at the original $(\omega,\omega)-$entry. Notice that $(1-t)(1+t+\cdots+t^{\omega-1})=1-t^{\omega}$. As a result,

\begin{equation}
|\left(R_{ji}(t)\right)|=(-1)^{\omega(e)-1}|\left(R_{ji}^{**}(t)\right)|, \label{for:Rijt replace t to tw}
\end{equation}
where $\left(R_{ji}^{**}(t)\right)$ is the minor of the original $\left(R_{ji}(t)\right)$ deleting the $1$st, $2$nd ,..., $(\omega-1)$th rows and the $1$st, $2$nd ,..., $(\omega-1)$th columns and replacing $t$ by $t^{\omega}$ in the original $\omega$th column.

If $\hat{\gamma}_{1}$ is a reverse arc of $P_{1}$, A similar argument shows that (\ref{for:Rijt replace t to tw}) still holds.

Do it inductively for path on $T$ corresponding to each weighted edge in $T^{\omega}$, we finally have

\begin{equation}
|\left(R_{ji}(t)\right)|=(-1)^{g-|T^{\omega}|}|\left(R_{ji}(t)\right)|, \nonumber
\end{equation}
where $g$ is the number of edges in $T$ and $|T^{\omega}|$ denotes the number of edges in $T^{\omega}$.
\end{proof}

\begin{example}
From the ribbon graph in Fig. \ref{fig:ribbon diagram ribbon graph}(2), we can contract $E_{4},E_{5}$ and $E_{6}$ to get a contracted ribbon graph. Then by formula (\ref{for:contracted Rijt}), the matrix $\left(\boldsymbol{R}_{ji}(t)\right)$ is

\begin{equation}
\boldsymbol{R}(t)=\left(\begin{array}[]{cccc}
1 & 1-t & 0 & 0 \\
0 & 1 & 1-t & 0 \\
t-1 & 0 & -t & 1-t^{3} \\
1-t & 0 & 0 & -1
\end{array}\right). \nonumber
\end{equation}

Thus the half Alexander polynomial is

\begin{equation}
A_{R}(t)=|\left(\boldsymbol{R}_{ji}(t)\right)|=t(1+t^{2}-3t^{3}+3t^{4}-t^{5}). \nonumber
\end{equation}
\end{example}

\subsection{Path-type formula, Theorem \ref{thm:algo3} and geography.}
We give a much simpler formula for half Alexander polynomial when the ribbon tree is a path. In fact, we can use finger moves to deform any ribbon to a new ribbon so that the ribbon tree is a path.

\subsubsection{Generality of path-type ribbon graph.}
\begin{prop}\label{prop:path}
Any ribbon knot has a ribbon such that $T$ is a path in the ribbon graph $(T,S)$.
\end{prop}

\begin{figure}[h]
\centering
\includegraphics[scale=0.25]{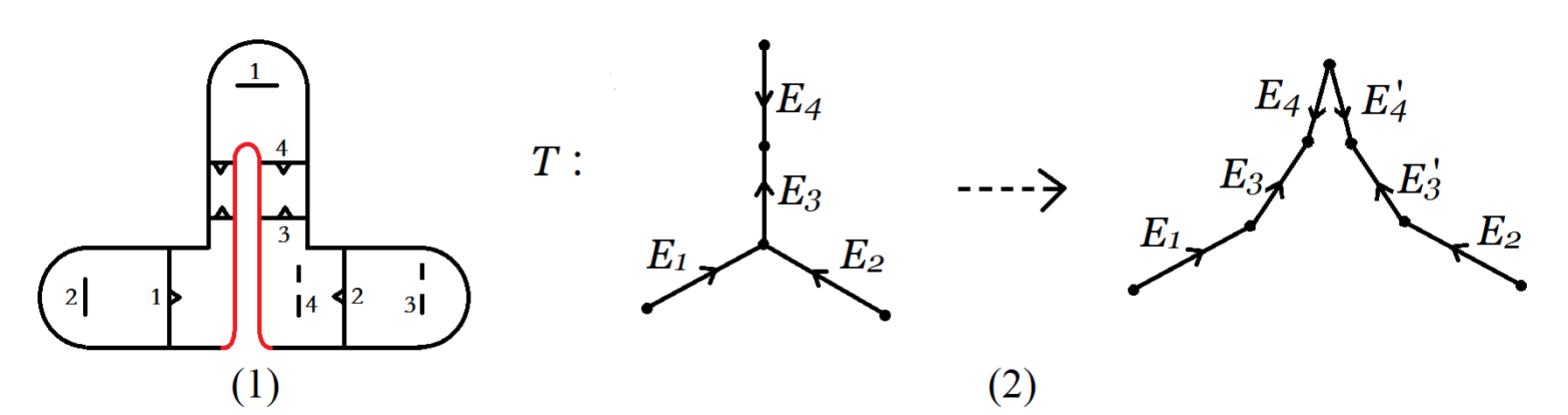}
\caption{Change ribbon graph into a path.}\label{fig:change ribbon graph into a path}
\end{figure}

\begin{proof}
Let $R$ be a ribbon in $S^{3}$, and $(D^{2},\cup_{i=1}^{g}\beta_{i},\cup_{i=1}^{g}\hat{\gamma}_{i})$ be its ribbon diagram. Note that a tree is a path if and only if it has only two pendant vertices. We use finger moves on $R$ to decrease pendant vertices.

Let $v_{i}$ be a pendant vertex of $T$. Let $P=E_{1}E_{2}...E_{k}$ be the maximal induced path in $T$ with one endpoint $v_{i}$ incident to $E_{1}$. Then the other endpoint denoted $v_{j}$, has degree at least 3 in $T$. Recall that each vertex of $P$ is a component of $D^{2}-\cup_{i=1}^{g}\gamma_{i}$. The boundary of the component corresponding to $v_{j}$ has an arc $a\subset r(\partial D^{2})$ not intersecting $\hat{\gamma}_{k}$, where $\hat{\gamma}_{k}$ corresponds the edge $E_{k}$. Push $a$ across $\alpha_{k},...,\alpha_{1}$ in succession by finger moves in $S^{3}$ to split each of $\hat{\gamma}_{k},...,\hat{\gamma}_{k}$ into two arcs, as shown in Fig. \ref{fig:change ribbon graph into a path}(1). The effect on $T$ is to split the edges $E_{k},...,E_{1}$ in succession until $v_{i}$ has degree 2, as shown in Fig. \ref{fig:change ribbon graph into a path}(2).

Do this inductively to pendant vertices until $T$ becomes a path.
\end{proof}

\subsubsection{Path-type formula.}

Let $R$ be a ribbon with ribbon graph $(T,S)$ where $T$ is a path. 
Assume $T=v_{1}v_{2}\cdots v_{n+1}$, and the directed edge $E_{i}$ is either $v_{i}v_{i+1}$ or $v_{i+1}v_{i}$ for $i=1,\cdots,n$. 
Read off a matrix $W_{(n+1)\times n}(t)=W^{\prime}(t)+W^{\prime\prime}(t)$ from $(T,S)$ as follows. Let
\begin{align}
\begin{cases}
W^{\prime}_{i,i}(t)=1,\,W^{\prime}_{i+1,i}(t)=-t,&if\,E_{i}=v_{i}v_{i+1},\\
W^{\prime}_{i,i}(t)=-t,\,W^{\prime}_{i+1,i}(t)=1,&if\,E_{i}=v_{i+1}v_{i},\\
W^{\prime\prime}_{k,i}(t)=t-1,&if\,S(E_{i})=v_{k}, \label{for:path type}
\end{cases}
\end{align}
and all the other entries be 0.

\begin{thm}\label{thm:algo3}
$A_{R}(t)=\pm|W^{*}(t)|$, where $W^{*}(t)$ is $W(t)$ deleting the $(n+1)th$ row.
\end{thm}

\begin{proof}
Let
\[
\Lambda_{n} = \begin{pmatrix}
1 & & & \\
1 & 1 & & \\
\vdots & \vdots & \ddots & \\
1 & 1 & \cdots & 1
\end{pmatrix}_{n \times n}, \qquad X = \begin{pmatrix}
sgnE_{1} & & & \\
& sgnE_{2} & & \\
& & \ddots & \\
& & & sgnE_{n}
\end{pmatrix},
\]
where $sgnE_{i}=+1/-1$ if and only if $E_{i}=v_{i}v_{i+1}/v_{i+1}v_{i}$. 
Then we have $R(t)=-X\Lambda_{n}W^{*}(t)$. 
One can verify this case by case depending on $sgnE_{i}$ and whether $E_{i}$ is forward in $P_{i}$. Therefore $A_{R}(t)=(-1)^{|E_{+1}|}|W^{*}(t)|$, where $|E_{+1}|$ is the number of edges with $sgn=+1$ in $T$.
\end{proof}

\begin{rem}
Theorem \ref{thm:algo3} has a more topological proof and $R(t)=-X\Lambda_{n}W^{*}(t)$ can be understood more topologically. We hint that one can refer to the proof of Proposition \ref{prop:relation with  Wirtinger} to find dual homological generators, but we will not explain the details here.
\end{rem}

\subsubsection{Geography of half Alexander polynomials.}

Given an integer-coefficient polynomial $f(t)$ with $f(1)=\pm 1$, Terasaka constructed in [15] a ribbon knot whose Alexander polynomial is $f(t)f(t^{-1})$. His arguments were complicated calculations using Wirtinger representation of the knot diagram. But one can verify that in fact his construction satisfies $A_{R}(t)\dot{=}f(t)$. We give a much more compact construction of ribbons such that $A_{R}(t)\dot{=}f(t)$ using Theorem \ref{thm:algo3}. But our calculation is still inspired by Terasaka's.

\begin{cor}\label{cor:geography}
For any $f(t)\in\mathbb{Z}[t]$ with $f(1)=\pm 1$, there is a ribbon diagram so that the half Alexander polynomial $A_{R}(t)=t^{|E_{-1}|}f(t)$, where $|E_{-1}|$ is the number of edges with $sgn=-1$ in the ribbon tree.
\end{cor}

\begin{proof}
First assume $f(1)=1$. Then $F(x)=(f(t)-1)/(t-1)$ is a polynomial. We may suppose $F(x)=\sum_{i=0}^{n}a_{i}t^{i}-\sum_{j=0}^{m}b_{j}t^{j}$, where $a_{i},b_{j}\geq 0$. We construct a ribbon graph $(T,S)$ where $T$ is a path as follows.

Let the two pendant vertices of $T$ be \textit{red vertex} and \textit{blue vertex}. In Fig. \ref{fig:geography ribbon graph}, $S$ maps each red/blue edge to the red/blue vertex. Start from red vertex. The first two edges are the blue edges as depicted in Fig. \ref{fig:geography ribbon graph}(1). Then identify the terminal vertices of the paths for $a_{0},...,a_{n}$ as shown in Fig. \ref{fig:geography ribbon graph}(2) in succession. If $n>m/n<m$, identify the terminal vertice of the path as shown in the upper/lower left of Fig. \ref{fig:geography ribbon graph}(3) and then identify the terminal vertices of the paths for $b_{m},...,b_{0}$ as shown in the upper/lower right of Fig. \ref{fig:geography ribbon graph}(3) in order. Ending at the blue vertex gives the path $T$.

\begin{figure}[h]
\centering
\includegraphics[width=0.8\textwidth]{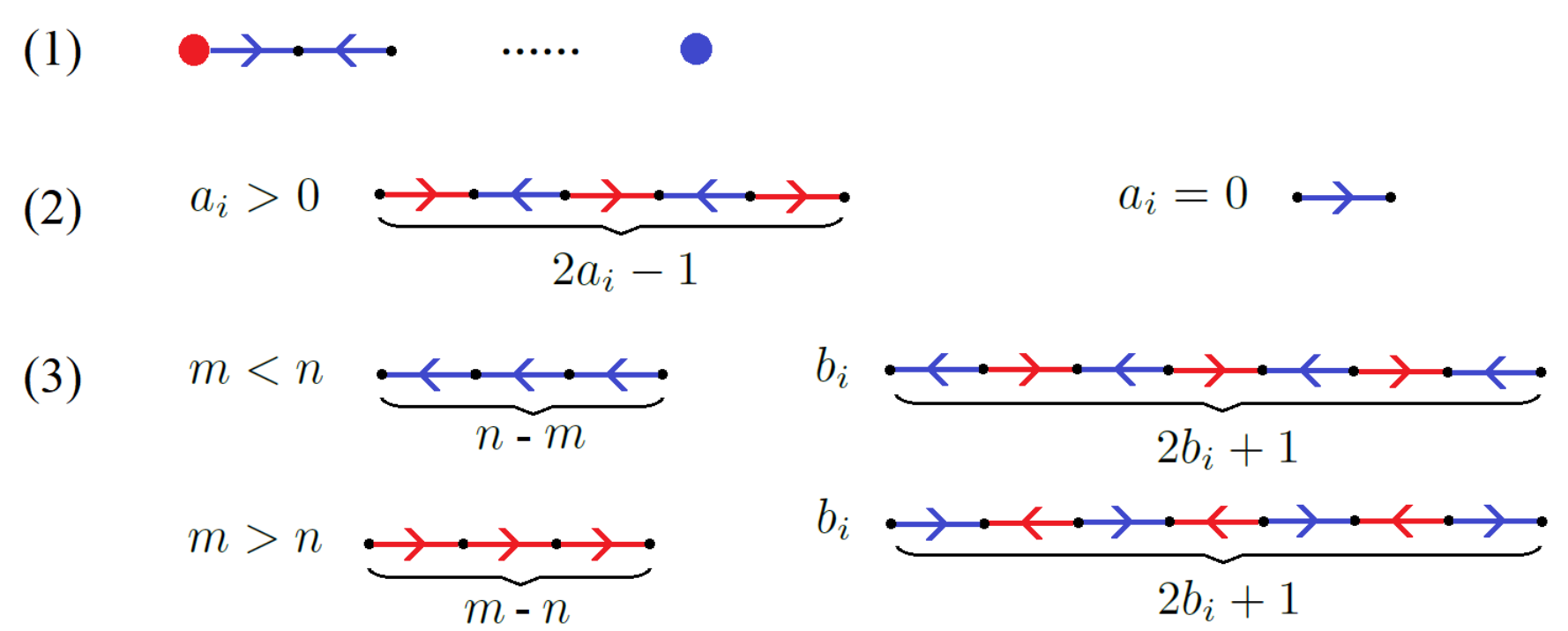}
\caption{The construction of the ribbon graph.}\label{fig:geography ribbon graph}
\end{figure}

To calculate $ A_R(t) $, observe the following two $( (2k+1) \times (2k+1) $ blocks.

\[
\begin{pmatrix}
t-1 & 0 & t-1 & \cdots & 0 & t-1 \\
1   &   &     &        &   &   \\
 -t & -t & & & & \\
& 1 & 1 & & & \\
& & -t & \ddots & & \\
& & & \ddots & -t & \\
& & & & 1 & 1
\end{pmatrix},
\begin{pmatrix}
0 & t-1 & 0 & \cdots & t-1 & 0 \\
-t & & & & & \\
1 & 1 & & & & \\
& -t & -t & & & \\
& & & 1 & \ddots & \\
& & & & \ddots & -t \\
& & & & & 1 & 1
\end{pmatrix}
\]
Multiplying the columns by $ 1, -1, 1, \cdots, -1, 1 $ in succession and adding together to the first column, the first column of the two blocks become $ ((k+1)(t-1), 1, 0, \cdots, 0)^T $ and $ (-k(t-1), 1, 0, \cdots, 0)^T $ respectively. The matrix $ W^*(t)$ defined in Theorem \ref{thm:algo3} is like

\[
\begin{pmatrix}
1 & 0 & t-1 & 0 & t-1 & t-1 & \cdots & 0 & t-1 & 0 \\
-t & -t & & & & & & & \\
& 1 & 1 & & & & & & \\
& & -t & -t & & & & & \\
& & & 1 & 1 & & & & \\
& & & & -t & 1 & & & \\
& & & & & -t & \ddots & & \\
& & & & & & \ddots & -t & \\
& & & & & & & 1 & 1 & -t
\end{pmatrix}
\]
Multiply the columns by $1,-1,1,-1,1,t,-t,...,t^{n},t^{n-1},...,t^{m},-t^{m},t^{m},...,1,-1,1$ in succession and add together to the first column, then the first column becomes

\[
\left(1+(t-1)(\sum_{i=0}^{n}a_{i}t^{i}-\sum_{j=0}^{m}b_{j}t^{j}),0,\cdots,0\right)^T.
\]

Thus $|W^{*}(t)|=(-t)^{|E-1|}\,f(t)$, where $|E_{-1}|$ is the number of edges with $sgn=-1$. Eventually $A_{R}(t)=|-X\Lambda_{n}W^{*}(t)|=t^{|E-1|}\,f(t)$.

If $f(1)=-1$, apply an inverse of R0-reduction to the ribbon diagram constructed above.
\end{proof}

\begin{figure}[htbp]
		  \centering
    \includegraphics[scale=0.21]{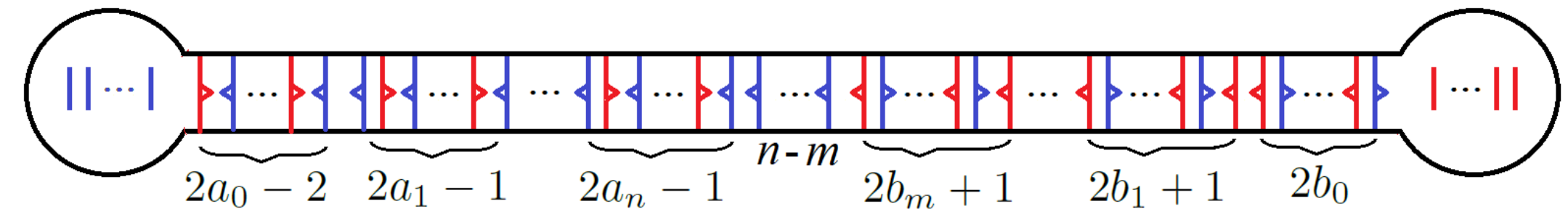}
    \caption{A ribbon with $A_R (t) = 1+(t-1)(\sum_{i=0}^{n}a_{i}t^{i}-\sum_{j=0}^{m}b_{j}t^{j}) $.}\label{fig:a ribbon geography}
	\label{fig:25}
\end{figure}

Actually, we can use R2- and R0- reductions to eliminate the first two and the last edge in $T$ in the above proof. See Fig. \ref{fig:a ribbon geography} for a typical example of corresponding ribbon.

\begin{rem}
We cannot get rid of $t^{|E-1|}$ in the statement of Corollary \ref{cor:geography}, as seen from the proof.
\end{rem}

\subsection{Topological meaning of half Alexander polynomial.}
It is basic to ask for a topological meaning of matrix $ R(t) $ and $ A_R(t) = |R(t)| $. We give a 3-dimensional explanation for them.

\begin{prop}\label{prop:nfolklore}
In the Alexander module of a ribbon knot $ K $ with ribbon $ R $, the submodule generated by $ e_1, ..., e_g $ has presentation matrix $ R(t^{-1})^T $, and thus $ A_R(t^{-1}) $ is a generator of the first elementary ideal, which is principal, of this submodule.
\end{prop}

\begin{proof}
Seen from (\ref{for:tAAT}), $-tR(t^{-1})^T$ is a presentation matrix of the submodule generated by $ e_1, ..., e_g$.
\end{proof}

Our explanation is essentially different from the folklore result that $ A_R(t)$ is the Alexander polynomial of the ribbon in $ B^4 $. In fact the latter come from the generators $ f_1, ..., f_g $ and $R(t)$ in (\ref{for:tAAT}).

\subsection{Half Alexander polynomial is a stronger invariant for ribbons.} \label{subsect:example}
We show, by a example with $\Delta_{K_1}(t) = \Delta_{K_2}(t)$ but $ A_{R_1}(t) \ne A_{R_2}(t)$, that half Alexander polynomial carries more information than Alexander polynomial for
ribbons even up to multiplication by $ \pm t ^{ \pm n} $ and changing $t$ by $t^{ - 1}$.

\begin{figure}[htbp]
		  \centering
    \includegraphics[scale=0.23]{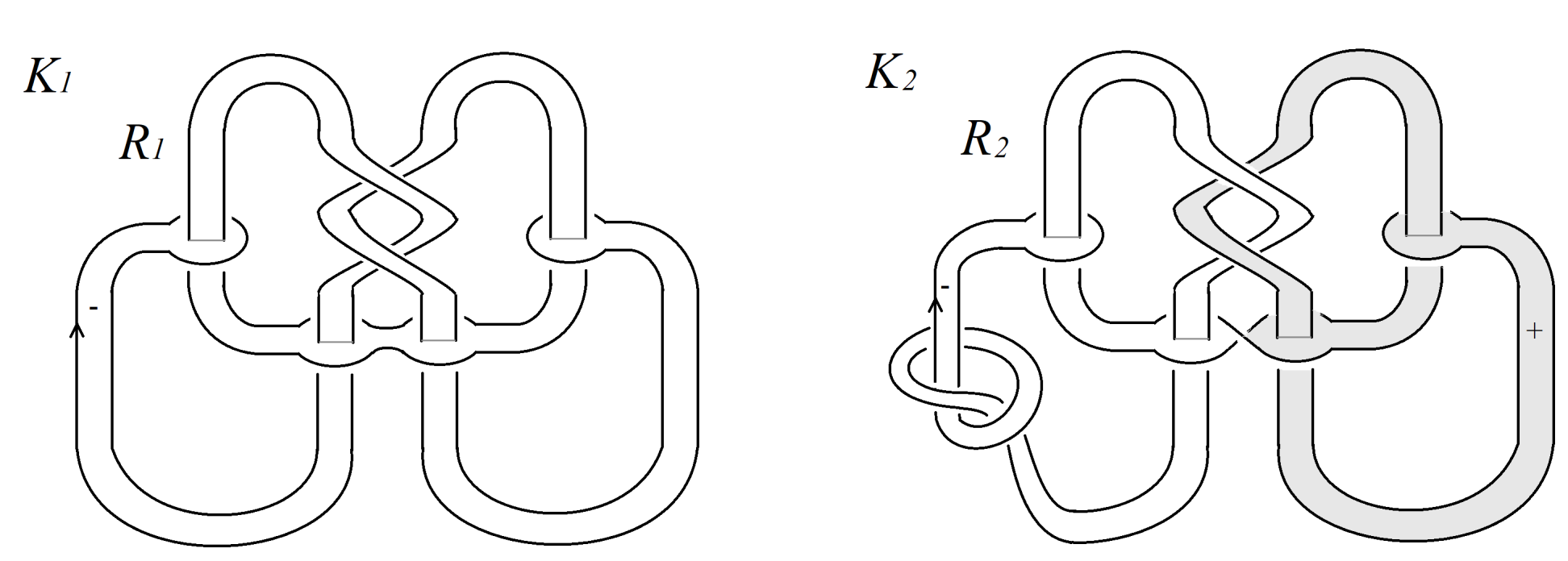}
    \caption{Ribbon knot $ K_1 $ with ribbon $ R_1 $ and ribbon knot $ K_2 $ with ribbon $ R_2 $.}
	\label{fig:half not equal Alexander polynomial equal}
\end{figure}

\begin{example}
There exist different ribbon knots with the same Alexander polynomial but having ribbons with different half Alexander polynomial. See Fig. \ref{fig:half not equal Alexander polynomial equal} for two ribbon knots and their ribbons. In fact, $A_{R_1}(t) = 1 - 4t + 4t^2$ and $A_{R_2}(t) = -2t + 5t^2 - 2t^3 $. But $ \Delta_{K_1}(t) = \Delta_{K_2}(t) = (2t^{-1} - 5 + 2t)^2 $. It can be verified that $ K_1 $ is hyperbolic, while $ K_2 $ is a satellite knot, so they are different knots.
\end{example}

\subsection{Half Alexander polynomial is not a ribbon knot invariant.}

The readers may guess a stronger result: half Alexander polynomial is an invariant for oriented ribbon knots up to multiplication by $\pm t^{\pm n}$. But this is incorrect. We find a ribbon knot from [3] and find two ribbons for it, whose half Alexander polynomials are different even up to multiplication by $\pm t^{\pm n}$ and changing $t$ by $t^{-1}$.

\begin{figure}[h]
\centering
\includegraphics[width=0.8\textwidth]{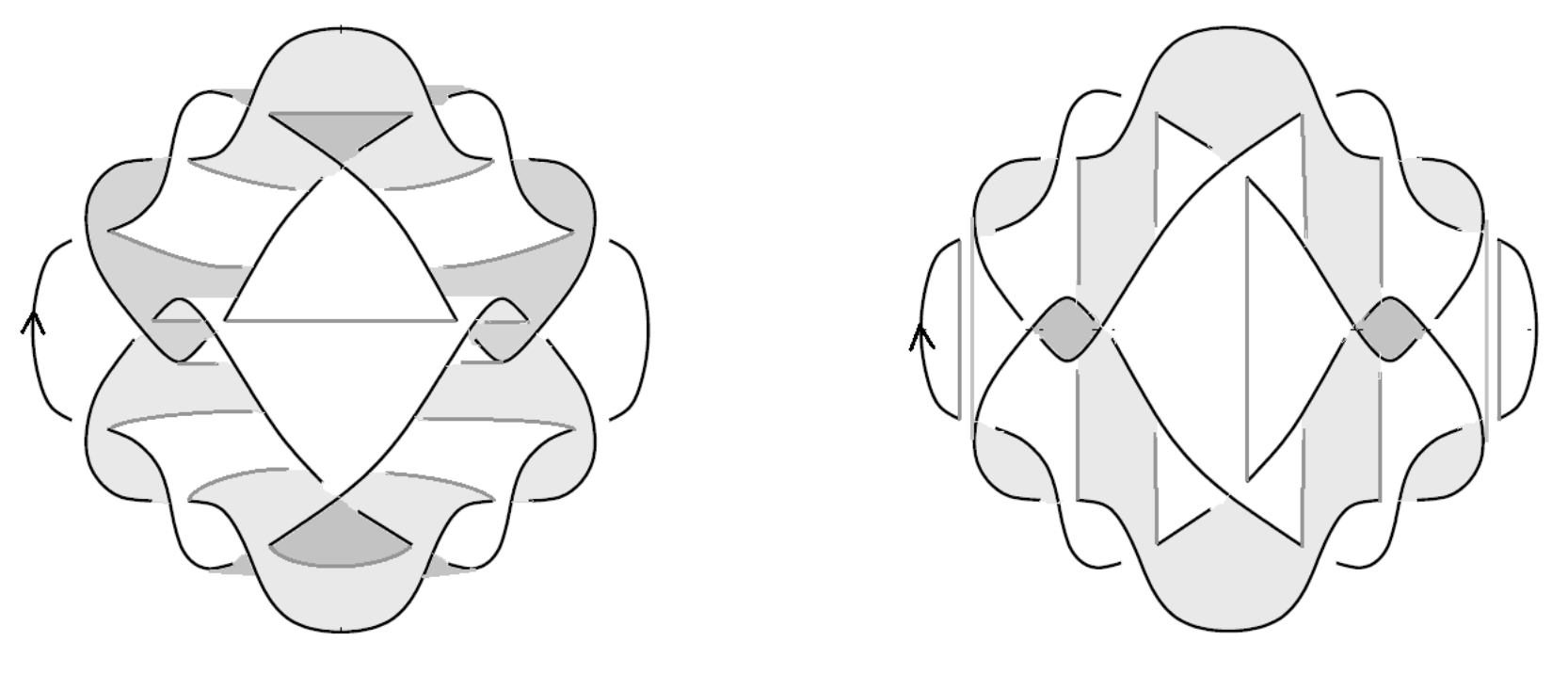}
\caption{Two different ribbons $R_{1}$ and $R_{2}$ for the knot $16{\rm n}524794$.}\label{fig:not ribon knot invariant}
\end{figure}

\begin{example}
For the knot $16{\rm n}524794$, consider the two ribbons in Fig. \ref{fig:not ribon knot invariant}. After calculation, we have $A_{R_{1}}(t)=-t^{3}(1-4t+4t^{2})$ while $A_{R_{2}}(t)=-2+5t-2t^{2}$.
\end{example}

This counterexample illustrates two facts. Firstly, by Lemma \ref{lem:reductions unchange} we obtain a geometric fact.

\begin{prop}
The two ribbons given in Fig. \ref{fig:not ribon knot invariant} are not related by ribbon 0-moves, ribbon 1-moves, finger moves or ribbon 3-moves, up to smooth isotopy of ribbons. $\square$
\end{prop}

Secondly, Fox and Milnor[4] proved that the Alexander polynomial of a slice knot has the form $f(t)f(t^{-1})$ for some integer-coefficient polynomial $f(t)$. This does not imply that $f(t)$ is a slice knot invariant because the Alexander polynomial may not factorize uniquely. But one might still expect that there exists a standard factorization of the Alexander polynomial such that $f(t)$ does become a slice knot invariant. To the best of our knowledge, no literature has explicitly denied this possibility. But our counterexample shows that this possibility is essentially nonexistent.

\section{A generalization of Theorem \ref{thm:algo2}. }

\begin{figure}
    \centering
    \includegraphics[scale=0.2]{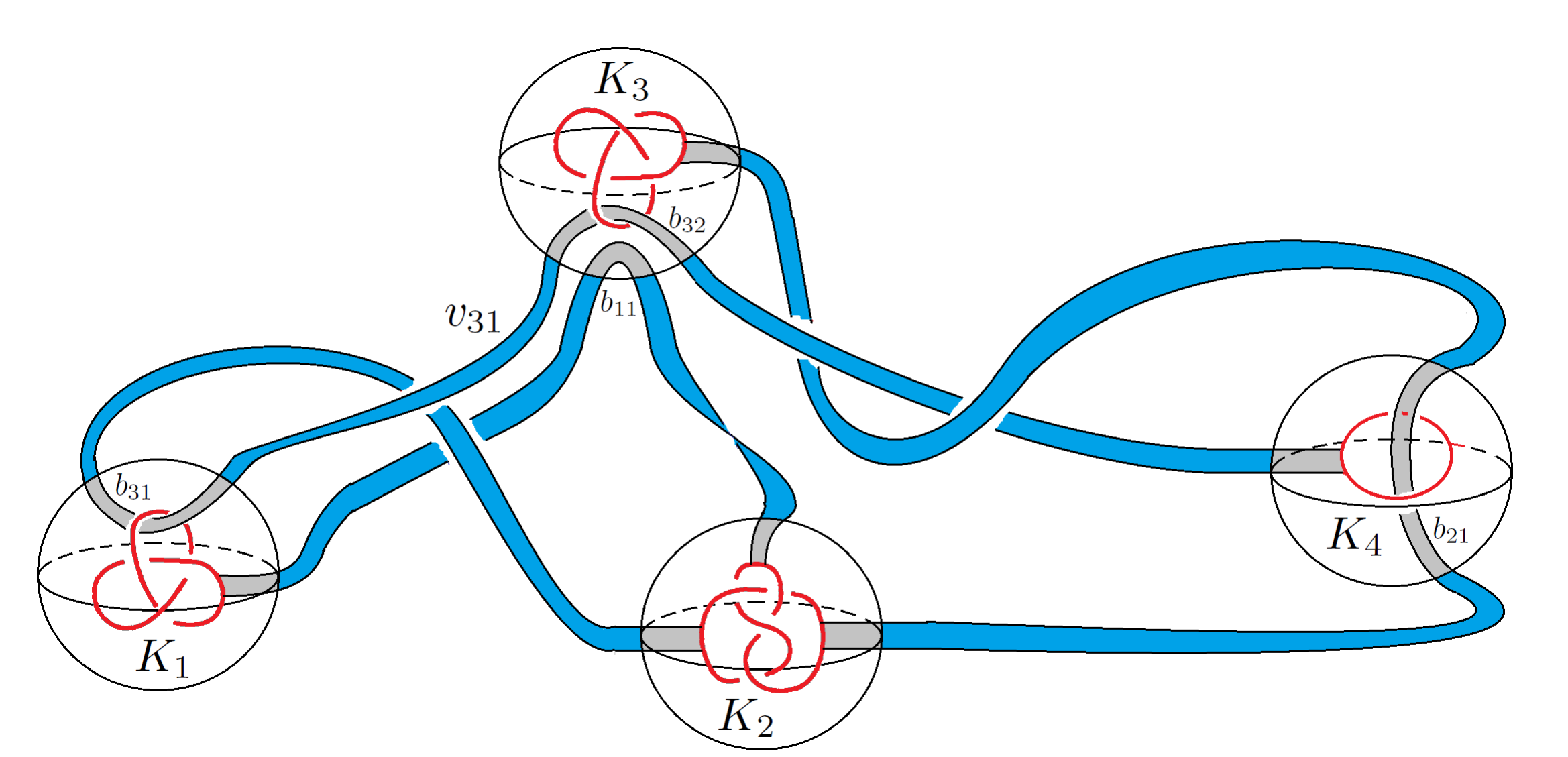}
    \caption{An example of $\cup_{i=1}^{n}K_{i}\bigcup\cup_{i=1}^{n-1} B_{i}$.}\label{fig:genelized knot example}
    \label{fig:33}
\end{figure}

Consider a split union of knots. To be specific, let $\Omega_{1},...,\Omega_{n}$ be disjoint balls in $S^{3}$ and $K_{i}$ be a knot (maybe unknot) in $\Omega_{i}$ for $i=1,...,n$. Let $B_{1},...,B_{n-1}$ be disjoint bands each connecting two knots with two arcs in the boundary so that $\cup_{i=1}^{n}K_{i}\bigcup\cup_{i=1}^{n-1}\partial B_{i}-int\left(\cup_{i=1}^{n}K_{i}\bigcap\cup_{i=1}^{n-1}\partial B_{i}\right)$ is a knot, denoted $K$. See Fig.\ref{fig:genelized knot example} for an example. We give a formula for Alexander polynomial of $K$.

Suppose $\cup_{i=1}^{n-1}\partial B_{i}$ intersects $\Omega_{1},...,\Omega_{n}$ transversely in arcs. Let $b_{i}$ be the core arc of $B_{i}$ connecting different knot components for $i=1,...,n-1$, oriented temporarily. Denote the components of $b_{i}\cap\cup_{i=1}^{n}\Omega_{i}$ in order along $b_{i}$ by $b_{i0},b_{i1},...b_{i,x_{i}}$. For each proper arc component $b_{ik}\subset\Omega_{j}$, there is a well-defined linking number $lk(b_{ik},K_{j})$. Denote the components of $b_{i}-\cup_{i=1}^{n}\Omega_{i}$ in order along $b_{i}$ by $v_{i1},...v_{i,x_{i}}$. We now define a \textit{contracted ribbon graph} $(T^{\omega},S)$, where $T^{\omega}$ is a weighted directed tree and $S$ is a map from its edges to its vertices, as follows.

\begin{itemize}
    \item \textit{Tree $T^{\omega}$}: The vertex set consists of vertex $V_{i}$ corresponding to $K_{i}$ for $i=1,...,n$ and vertices $v_{i1},...v_{i,x_{i}}$ for each $i=1,...,n-1$. Each $b_{ik}$ gives an edge $E_{ik}$ incident to the two vertices containing the endpoints of $b_{ik}$. Assume $b_{ik}\subset\Omega_{j}$, then the direction of $E_{ik}$ is from $v_{i,k-1}$ to $v_{ik}$ if $lk(b_{ik},K_{j})>0$ and opposite if otherwise. The weight of $E_{ik}$, denoted $\omega_{ik}$, is $|lk(b_{ik},K_{j})|$. Contract each edge $E_{i0}$ and each edge $E_{i,x_{i}}$.
    \item \textit{Singularity map $S$}: $S(E_{ik})=V_{j}$ if and only if $b_{ik}\subset\Omega_{j}$.
\end{itemize}

\begin{figure}
    \centering
    \includegraphics[scale=0.23]{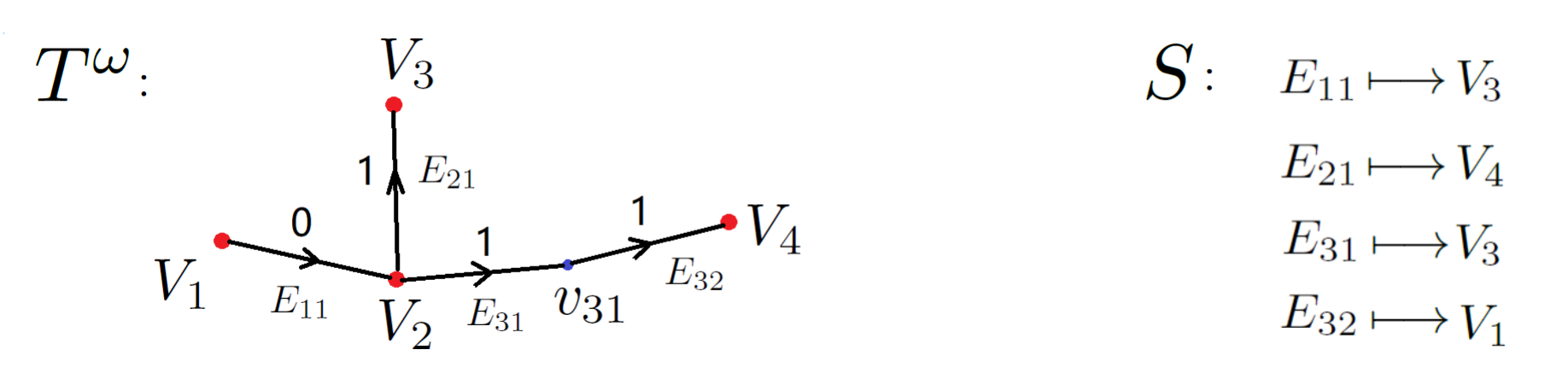}
    \caption{Contracted ribbon graph.}\label{fig:contracted ribbon graph genelized knot example}
    \label{fig:34}
\end{figure}

For example, Fig. \ref{fig:contracted ribbon graph genelized knot example} is the contracted ribbon graph for Fig.\ref{fig:genelized knot example}. Then we define a matrix $(\boldsymbol{R}_{ij}(t))$ from $(T^{\omega},S)$ according to (\ref{for:contracted Rijt}) and the paragraph before it.

\begin{thm}\label{thm:generalize}
With the above notations, the Conway-normalized Alexander polynomial of $K$ is
\[
\Delta_{K}(t)=|\boldsymbol{R}_{ij}(t)||\boldsymbol{R}_{ij}(t^{-1})|\Delta_{K_{1}}(t)\cdots\Delta_{K_{n}}(t).
\]
Especially, Alexander polynomials of $K$ is determined by the Alexander polynomials of $K_{1},...,K_{n}$ and the contracted ribbon graph.
\end{thm}

Rather than give a detailed proof we outline it. We will prove a more general result in detail including this theorem as a special case in \cite{Bai2023}.

\begin{proof}
First, let $F_{i}\subset\Omega_{i}$ be a Seifert surface for $K_{i}$ so that $\cup_{i=1}^{n}F_{i}\bigcup\cup_{i=1}^{n-1}B_{i}$ is a immersed surface with only ribbon singularities. Suppose the subband of $B_{i}$ corresponding to $b_{ik}\subset\Omega_{j}$ has ribbon arcs $\beta_{ik1},...\beta_{ik,y_{ik}}$. Desingularize the immersed surface into a Seifert surface for $K$, denoted $F_{K}$, as in Subsection \ref{sss:From ribbon to Seifert surface}.

Second, let $h_{i1},...,h_{i,z(i)}$ be a basis of $H_{1}(F_{i})$. Let $e_{ik,l}$ be a small circle encircling $\beta_{ik,l}$. Choose $f_{ik,l}$ to be the path on the immersed surface and $f_{ik,l}$ be the simple closed curve on $F_{K}$ as in Subsection \ref{sss:Homological basis on Fg}.

Then all of $e_{\cdots}$ , all of $f_{\cdots}$ and $h_{11},...,h_{1,z(1)},...,h_{n1},...,h_{n,z(n)}$ form a basis of $H_{1}(F_{K})$. The Seifert matrix of $F_{K}$ has the form
\begin{align*}
A = \left(
\begin{array}{ccccc}
O & P &     & 0 &\\
Q & L &     & \ast &\\
  &   &      A_1 & & \\
O & \ast & &\ddots & \\
  &   &      & & A_n
\end{array}
\right),
\end{align*}
where $A_{i}$ is a Seifert matrix for $K_{i}$.

Therefore
\[
\Delta(t)=|t^{\frac{1}{2}}A-t^{-\frac{1}{2}}A^{T}|=|tP-Q^{T}||t^{-1}P^{T}-Q|\Delta_{K_{1}}(t)\cdots\Delta_{K_{n}}(t).
\]

We can obtain $P,Q$ exactly as in Subsection \ref{fig:homological basis locally represented} to get formula (\ref{for:matrix R(t)}).

Finally, preceding exactly as in Subsection \ref{subsect:Contracted formula}, we have $|\boldsymbol{R}_{ij}(t)|=|tP-Q^{T}|$.
\end{proof}

\section{Alexander polynomials of general knots and virtual knots}

\subsection{Inspiration: Similarity between ribbon graph and Gauss diagram.}\label{subsect:inspiration}

Let $K$ be an oriented knot in $S^{3}$ with a diagram $\mathbb{K}$. Let $-\bar{K}$ denote the reflection of $K$ with the reversed orientation. Then the connected sum $K\#-\bar{K}$ is a ribbon knot, and we get a canonical ribbon from the diagram $\mathbb{K}\#-\bar{\mathbb{K}}$.

\begin{figure}[h]
\centering
\includegraphics[width=0.8\textwidth]{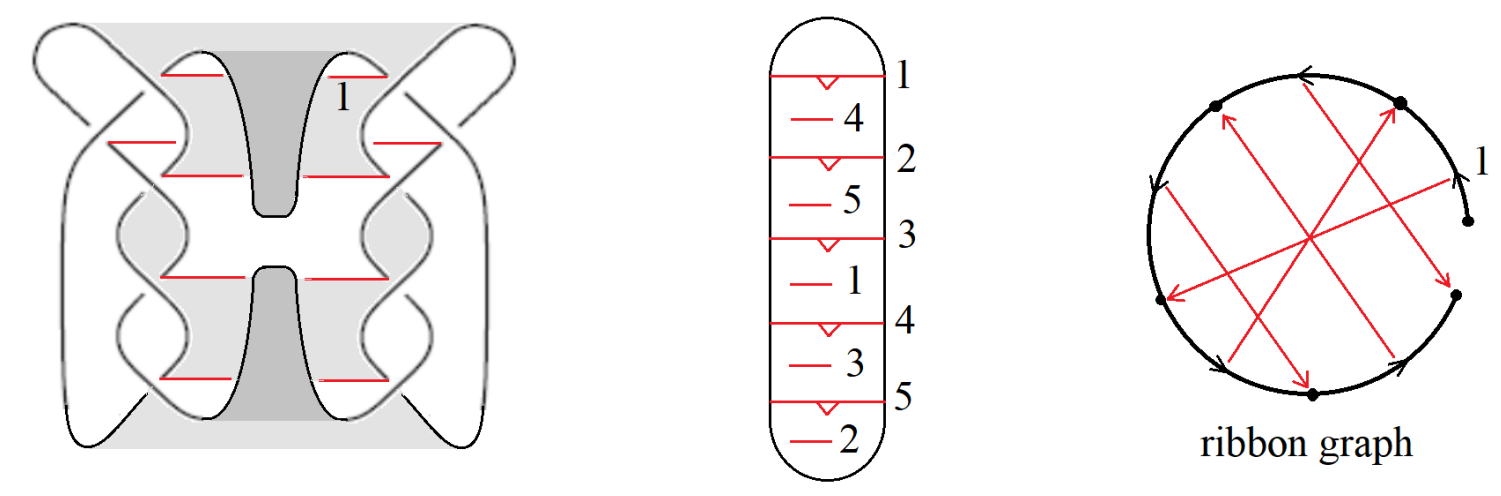}
\caption{Gauss diagram and ribbon graph.}\label{fig:Gauss diagram and ribbon graph}
\end{figure}

As illustrated in Fig. \ref{fig:Gauss diagram and ribbon graph} and \ref{fig:Gauss diagram}, from the Gauss diagram of $\mathbb{K}$ we can always get exactly a ribbon graph as follows.

\begin{enumerate}
\item Replace the positive/negative sign at the start of each arrow by a direction on the circle consistent/inconsistent with the orientation of $K$;
\item Merge the contiguous tips of arrows into a vertex;
\item Erase a small arc where the connected sum is made.
\end{enumerate}

We know that

\[
\Delta_{K\#-\bar{K}}(t)=\Delta_{K}(t)\Delta_{-\bar{K}}(t)=\Delta_{K}(t)\Delta_{K}(t)=\Delta_{K}(t)\Delta_{K}(t^{-1}).
\]

On the other hand, $\Delta_{K\#-\bar{K}}(t)=A_{R}(t)A_{R}(t^{-1})$, where $A_{R}(t)$ is the half Alexander polynomial of the ribbon. One may wonder whether $A_{R}(t)=\pm t^{\pm k}\Delta_{K}(t)$? We prove this is true. This leads us to discover a new method for computing Alexander polynomials of general knots.

\subsection{New formula for classical knots.} \label{subsect: classical knots}

\begin{figure}[h]
\centering
\includegraphics[width=0.6\textwidth]{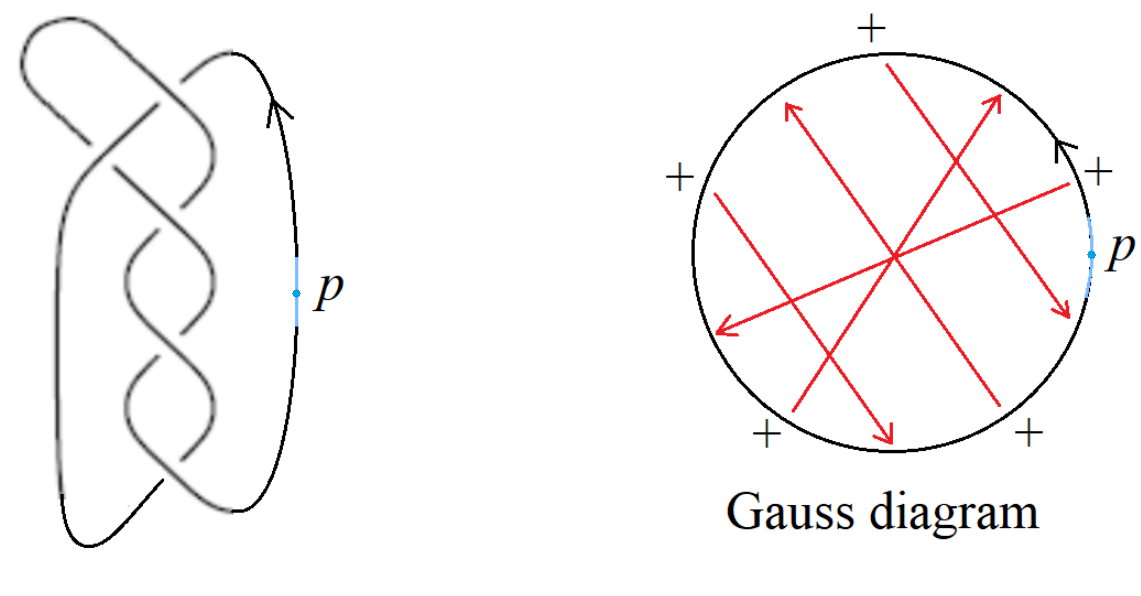}
\caption{Gauss diagram.}\label{fig:Gauss diagram}
\end{figure}

Let $\mathbb{K}$ be a diagram for a knot $K$ in $S^{3}$, with $n$ crossings $c_{1},...,c_{n}$. Let $\mathbb{G}$ be the Gauss diagram of $\mathbb{K}$. Denote the circle of $\mathbb{G}$ by $C$. Each crossing $c_{i}$ corresponds to an arrow, still denoted $c_{i}$. Choose a base point $p$ on $C$ avoiding the starts and tips of all the arrows. We adopt the convention that the start of the arrow corresponds to the lower strand and the tip of the arrow corresponds to the upper strand. See Fig. \ref{fig:Gauss diagram} for an example. We give formula for computing Alexander polynomial of $K$ from Gauss diagram $\mathbb{G}$. 

Denote the start of $c_{i}$ by $c_{i0}$. For each arrow $c_{i}$, let $P_{i}$ be the unique arc in $C-\{p\}$ from the start of $c_{i}$ to the tip of $c_{i}$, and define the sign of $P_{i}$, denoted $sgnP_{i}$, to be $+1$ if $P_{i}$ is counterclockwise and $-1$ otherwise.

Read off a matrix $\rho_{n\times n}$ as follows.
\begin{align}
    \rho_{ii}  &   = \frac{1}{2}\cdot sgnc_{i}\cdot sgnP_{i}; \nonumber\\ 
    \rho_{ji} & = \begin{cases}
                      sgnc_{j}\cdot sgnP_{i},&if~c_{j0}\in P_{i};\\
                       0,     &if~c_{j0}\notin P_{i}.
                   \end{cases}\forall j\neq i.
\end{align}
Then let $ R(t):=(t-1)\rho-\frac{1}{2}(t+1)I $.

\begin{thm}\label{thm:classical}
For the Alexander module of $K$, $R(t)$ is a presentation matrix.
Especially, the Alexander polynomial $\Delta_K (t)\dot{=} |R(t)|$.
\end{thm}

In fact, our formula is
\begin{align}
R_{ii}(t) &= \begin{cases}
-1, & \text{if } sgnP_{i} = sgnc_{i}; \\
-t, & \text{if } sgnP_{i} = -sgnc_{i}.\\ \nonumber
\end{cases}\\
R_{ji}(t) &= \begin{cases}
t-1, & \text{if } c_{j0} \in P_{i} \text{ and } sgnP_{i} = sgnc_{j}; \\
1-t, & \text{if } c_{j0} \in P_{i} \text{ and } sgnP_{i} = -sgnc_{j}; \\
0, & \text{if } c_{j0} \notin P_{i}.
\end{cases} \quad \forall j \neq i. \label{for:R(t) classicla knots}
\end{align}

\begin{rem}
According to the relationship between Gauss diagram and ribbon graph in Subsection \ref{subsect:inspiration}, R0-reduction applies to simplify the matrix $R(t)$. So to simplify the calculation, we can choose base point $p$ to be adjacent to a string of consecutive undercrossings.
\end{rem}

\begin{example}
For the knot $5_{2}$ in Fig. \ref{fig:Gauss diagram}, we have
\[
\rho=\left(\begin{array}[]{ccccc}
1/2&0&0&0&0\\
1&1/2&0&-1&0\\
1&1&1/2&-1&-1\\
0&1&1&-1/2&-1\\
0&1&0&0&-1/2
\end{array}\right),
\]
thus $\Delta(t){\dot{=}}|R(t)| = |(t-1)\rho-\frac{1}{2}(t+1)I|=2t^{-1}-3+2t$.
\end{example}

We can improve the matrix $R(t)$ to give a formula for the Conway-normalized Alexander polynomial. With the notations in the first two paragraphs in this subsection, define $\mathcal{R}(t)$ as follows.
\begin{align}
    \mathcal{R}_{ii}(t) &  =1 .\nonumber\\
    \mathcal{R}_{ji}(t)&=\begin{cases}
1-t^{sgnP_i sgnc_i}, & \text{if } c_{j0} \in P_i \text{ and } sgn c_i = sgn c_j; \\
t^{sgnP_i sgnc_i}-1, & \text{if } c_{j0} \in P_i \text{ and } sgn c_i = -sgn c_j; \quad \forall j \neq i.\\ \label{for:normal R(t) classical}
0,  & \text{if } c_{j0} \notin P_i.
\end{cases}
\end{align}
\begin{cor}\label{cor:nomal Alexander polynomial classical knot}
    The Conway-normalized Alexander polynomial $\Delta_K (t) = | \mathcal{R}(t) |$.
\end{cor}

\begin{proof}
Recall that $|R(t)|\doteq\Delta_K(t)$. As $|R(1)|=(-1)^n$, we have $|R(t)|=(-1)^n t^k \Delta(t)$ for some $k \in \mathbb{N}$. We claim that $k$ is the number of diagonal elements $-t$.

In fact, we know that $\Delta(1)=1$ and $\Delta'(1)=0$, so $|R(t)'| \big |_{t=1}=(-1)^n [k t^{k-1} \Delta(t) + t^k \Delta'(t)] \big |_{t=1}=(-1)^n k$. 
On the other hand, taking the derivative of $R(t)$ column by column and letting $t=1$, by (\ref{for:R(t) classicla knots}), we obtain $|R(t)'| \big |_{t=1}$ is $(-1)^n$ multiplying the number of diagonal elements $-t$.
Now dividing each column of $R(t)$ by the corresponding diagonal element $-1$ or $-t$, we obtain a matrix $\mathcal{R}(t)$.
\end{proof}

\begin{proof}[Topological and direct Proof of Theorem \ref{thm:classical}]
Choose a ball $B^{3}$ in $S^{3}$ so that $K-intB^{3}$ is an unknotted arc in $S^{3}-intB^{3}$. Set $K^{*}=K\cap B^{3}$. Then $B^{3}-N(K^{*})\cong S^{3}-N(K)$.

Assume $K^{*}$ is contained in a collar $N_{B}$ of $B^{3}$ and $\mathbb{K}^{*}$ is a regular projection of $K^{*}$ on $S^{2}=\partial B^{3}$. Orient $K^{*}$. Smooth the crossings of $\mathbb{K}^{*}$ in the orientation-preserving way to get Seifert circles and a Seifert arc on $S^{2}$, as shown in Fig. \ref{fig:The construction of surface F}(2).

Now we construct a surface $F$ in $N_{B}$ having the shape as drawn in Fig. \ref{fig:The construction of surface F}(3) near each crossing $c_{i}$ and Fig. \ref{fig:The construction of surface F}(4) near the two points $\partial K^{*}$, whose boundary is the union of $K^{*}$, the Seifert circles and Seifert arc.

Let $X_{\infty}$ be the infinite cyclic cover of $B^{3} - N(K^{*}) \cong S^{3} - N(K)$. Let $Y$ be the space $\left(B^{3} - N(K^{*})\right)$-cut-along-$F$, that is, $B^{3} - N(K^{*}) - F$ with two copies $F_{+}$ and $F_{-}$ replacing the removed $F$. We claim that $Y$ is a fundamental domain for the covering space $X_{\infty}$.

\begin{figure}
    \centering
    \includegraphics[scale=0.25]{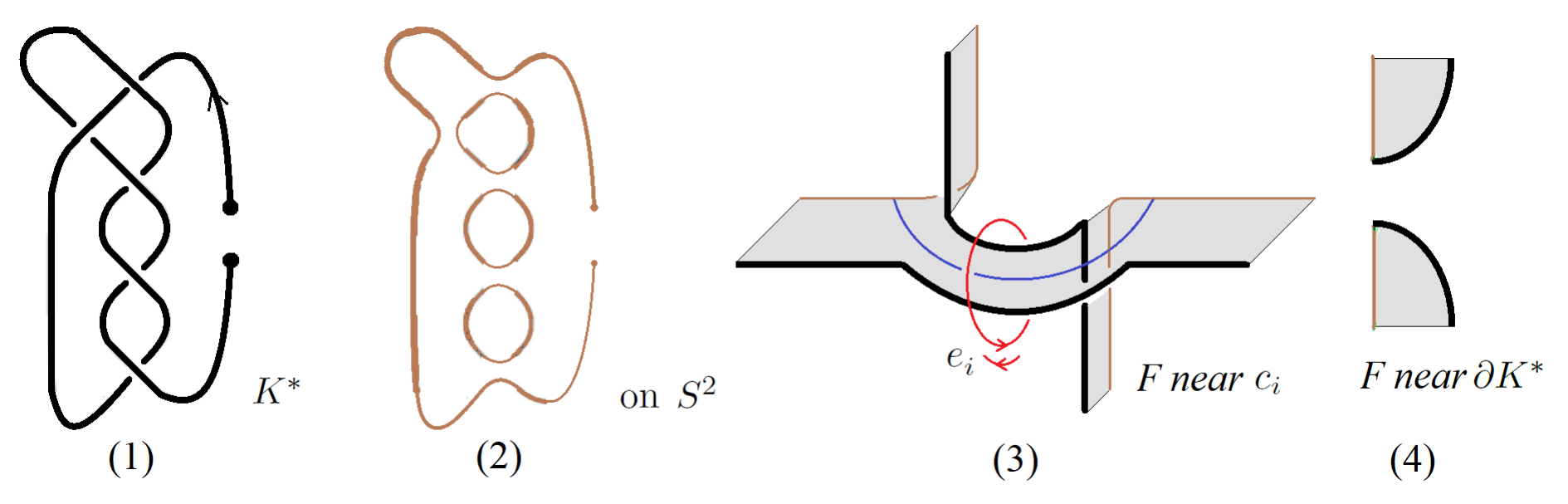}
    \caption{The construction of surface $F$.}\label{fig:The construction of surface F}
    \label{fig:29}
\end{figure}

In fact, the boundary of $B^{3} - N(K^{*})$ is a torus and the Seifert circles on $S^{2}$ are trivial circles on this torus. Cap these circles by the disks they bound on $S^{2}$ one by one from innermost ones and perturb it into the interior of $B^{3} - N(K^{*})$, then we get a Seifert surface $\bar{F}$ for $K$. It is well-known that $\left(B^{3} - N(K^{*})\right)$-cut-along-$\bar{F}$ is a fundamental domain of $X_{\infty}$. We identify the difference between it and $\left(B^{3} - N(K^{*})\right)$-cut-along-$F$. The latter subtracting the former is the union of solid neighbourhoods of the Seifert disks on the boundary torus, as 1-handles connecting the boundary torus to $\bar{F}$, while the former subtracting the latter is the same solids, but as 2-handles. Thus in $X_{\infty}$, for any $i \in \mathbb{Z}$, cut off these 2-handles in the $i$th copy of $\left(B^{3} - N(K^{*})\right)$-cut-along-$\bar{F}$ and glue them to the corresponding adjacent copy as 1-handles, then we get copies of $\left(B^{3} - N(K^{*})\right)$-cut-along-$F$. See Fig. \ref{fig:Seifert surface vs cut along surface}.

\begin{figure}
    \centering
    \includegraphics[scale=0.25]{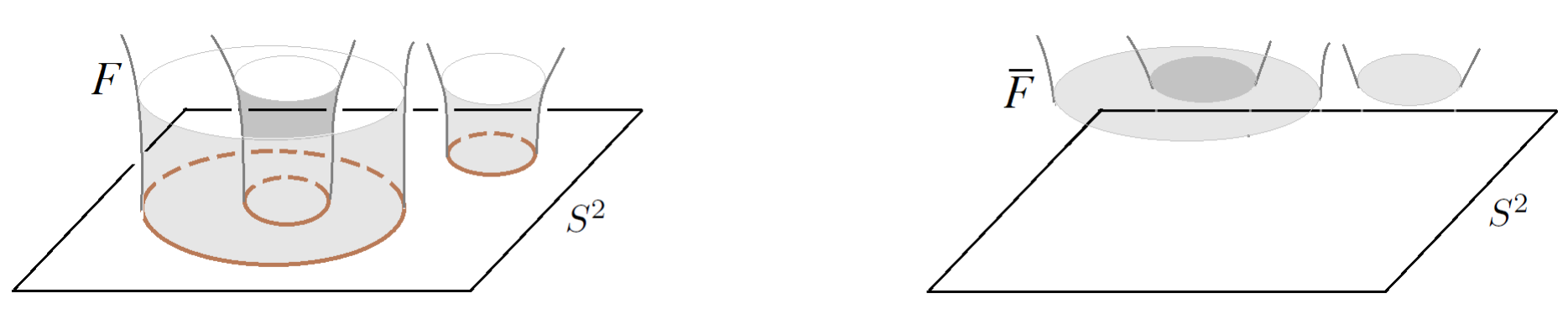}
    \caption{$\left(B^{3} - N(K^{*})\right)$-cut-along-$F$ vs cut-along $\bar{F}$.}\label{fig:Seifert surface vs cut along surface}
    \label{fig:30}
\end{figure}

For each crossing $c_{i}$, let $e_{i}$ be a simple closed curve in $B^{3} - F$ encircling the band as shown in Fig. \ref{fig:The construction of surface F}(3), oriented to encircle the upper strand of $K^{*}$ with a right-hand screw. Then $e_{1},...,e_{n}$ form a basis of $H_{1}(B^{3} - F)$.

On the other hand, if we cut along all the blue arcs depicted in Fig. \ref{fig:The construction of surface F}(3) for the crossings, then $F$ becomes a disk. So we may choose simple closed curves $f_{1},...,f_{n}$ on $F$ dual to the blue arcs to form a basis of $H_{1}(F)$. Specifically, $f_{i}$ starts from the lower strand at $c_{i}$, ends at the upper strand at $c_{i}$, and intersects the blue arc at $c_{i}$ in a single point. See Fig. \ref{fig:the basis cut long blue arcs}. Define the sign of $f_{i}$, denoted $sgnf_{i}$, to be $+1$ if the orientation of $f_{i}$ is similar to that of $K^{*}$ and $-1$ otherwise.

\begin{figure}
    \centering
    \includegraphics[scale=0.25]{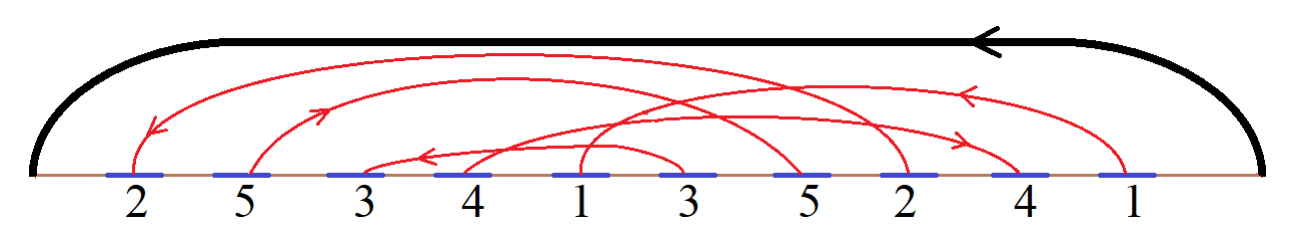}
    \caption{The basis $f_{1},...,f_{n}$ drawn on $F$ cut along the blue arcs.}\label{fig:the basis cut long blue arcs}
    \label{fig:31}
\end{figure}

Notice that $F$ is orientable, with the orientation compatible with that of $K^{*}$. Let $f_{i}^{+}/f_{i}^{-}$ be the push-off of $f_{i}$ in the positive/negative normal direction of $F$ into $B^{3} - F$. Then by tracing along $f_{i}$ in Fig. \ref{fig:The construction of surface F}(3), we see that, in $H_{1}(B^{3} - F)$,
\begin{align}
 [f_{i}^{+}] = \Sigma_{j=1}^{n} P_{ji}[e_{j}], \qquad [f_{i}^{-}] = \Sigma_{j=1}^{n} Q_{ij}[e_{j}]    \label{for:f+- P Q}
\end{align}
where
\[
P_{ii} = \begin{cases}
0, & \text{if } sgnf_{i} = sgnc_{i}; \\
-1, & \text{if } sgnf_{i} = -sgnc_{i}.
\end{cases}
\]
\[
Q_{ii} = \begin{cases}
+1, & \text{if } sgnf_{i} = sgnc_{i}; \\
0, & \text{if } sgnf_{i} = -sgnc_{i}.
\end{cases}
\]
\[
\forall j \neq i, \quad P_{ji} = Q_{ij} = \begin{cases}
sgnc_{j} \cdot sgnf_{i}, & \text{if } f_{i} \text{ goes through the lower strand at } c_{j}; \\
0, & \text{otherwise}.
\end{cases}
\]

The deck transformation on $X_{\infty}$ corresponding to the meridian of $K$ induces an automorphism $t$ on $H_{1}(X_{\infty})$. By (\ref{for:f+- P Q}), it is standard to prove that $R(t) = tP - Q^{T}$ is a presentation matrix for $\mathbb{Z}[t,t^{-1}]$-module $H_{1}(X_{\infty})$, that is, the Alexander module, where $P = (P_{ji})_{n \times n}$ and $Q = (Q_{ji})_{n \times n}$. Using notation in terms of Gauss diagram, we have (\ref{for:R(t) classicla knots}).
\end{proof}

\subsection{Relationship known formulae}
Although we have proved Theorem \ref{thm:classical} independently, it is still thought-provoking to make clear the relation between our formula and as many known formulas as possible.

\subsubsection{Canonical Seifert surface method.}
Although $F$ constructed in the above proof is similar to the canonical Seifert surface, it is difficult to find a general conversion relationship between our method and the Seifert surface method. The two matrices may be quite different, and the size of Seifert matrix for the canonical Seifert surface is less than $n$. Furthermore, our construction of $F$ is valid for virtual knots in general, while a major difficulty in defining Alexander polynomials for virtual knots is that they cannot bound Seifert surfaces.

\subsubsection{Wirtinger presentation formula.}
We can make clear the relation between the two presentation matrices $R(t) = tP - Q^{T}$ and $W = \Phi - t\Psi$ in Subsection \ref{subsect:virtual Wirtinger} for the Alexander module.

There is no direct linear translation between $R(t)$ and $W$. We need to enlarge $W$ depending on the choice of the base point $p$. Assume the base point $p$ is in the over-passing arc between $c_{n}$ and $c_{1}$, and denote the generators for the initial and last over-passing arc of $\mathbb{K} - p$ by $x_{0,1}$ and $x_{n,n+1}$ respectively, which are in fact equal.

Recall that Fig. \ref{fig:FoxofWirtinger} gives $W = \Phi - t\Psi$. Now let $\tilde{W} = \dot{\Phi} - t\tilde{\Psi}$ be the matrix representing the relations in Fig. \ref{fig:FoxofWirtinger} for the long knot $\mathbb{K} - p$, namely,
\[
\dot{\Phi}_{i+\frac{1}{2}(1-sgnc_{i}),i } = -1, \quad \dot{\Phi}_{k(i)+1,i} = 1;
\]
\[
\dot{\Psi}_{i+\frac{1}{2}(1+sgnc_{i}),i } = -1, \quad \dot{\Psi}_{k(i)+1,i} = 1,
\]
and other entries are $0$. 
Then $\tilde{W} = \tilde{\Phi} - t\tilde{\Psi}$ is an enlarged and equivalent presentation matrix of $W = \Phi - t\Psi$, where $\tilde{\Psi} = \left(\eta \quad \dot{\Phi}_{(n+1) \times n}\right)$, $\tilde{\Psi} = \left(O_{(n+1) \times 1} \quad \dot{\Psi}_{(n+1) \times n}\right)$ and $\eta = (-1,0 \cdots 0,1)^{T}$.

Set
\[
\Lambda_{n} = \begin{pmatrix}
1 & & & \\
1 & 1 & & \\
\vdots & \vdots & \ddots & \\
1 & 1 & \cdots & 1
\end{pmatrix}_{n \times n}, \qquad X = \begin{pmatrix}
sgnc_{1} & & & \\
& sgnc_{2} & & \\
& & \ddots & \\
& & & sgnc_{n}
\end{pmatrix}.
\]
Our general transformation formula is 
\begin{prop}\label{prop:relation with  Wirtinger}
$R(t) = -X\Lambda_{n}\dot{W}^{*}$, where $\dot{W}^{*} = \dot{\Phi}^{*} - t\dot{\Psi}^{*}$ is $\dot{W} = \dot{\Phi} - t\dot{\Psi}$ deleting the $(n+1)$th row.
\end{prop}

\begin{rem}
It is always incorrect that $R(t) = -X\Lambda_{n}W$, as $|W| = 0$.
\end{rem}

This can be verified directly. But we prefer to derive it in a more comprehensible way.

\begin{proof}[Proof of Proposition \ref{prop:relation with  Wirtinger}]
Choose generators $x_{i-1,i} - x_{n,n+1}$, $i = 1,...,n$ for Alexander module, then $\tilde{W}$ deleting the last row is a presentation matrix for Alexander module. We further have the following subtle observation.

\textbf{Claim 1:} For Alexander module, $\dot{W}^{*}$ is a presentation matrix.

In fact, the first relation $x_{0,1} - x_{n,n+1} = 0$ is a linear combination of the relations from the second to the $n+1$th columns multiplying some suitable $t^{k_{i}}$'s. The proof is similar to the well-known proof that in $W = \Phi - t\Psi$, one column is a linear combination of the other columns multiplied with suitable $t^{k_{i}}$'s.

Now Set element $e_{i} = x_{k(i),k(i)+1} - x_{i-sgnc_{i},i}$ for $i = 1,...,n$, which in fact is the $e_{i} $ in Fig. \ref{fig:The construction of surface F}(3).

\textbf{Claim 2:} For the Alexander module, $e_{1}, \cdots, e_{n}$ form generators.

In fact, introduce equivalent generators $\delta_{i} = x_{i,i+1} - x_{i-1,i}$ for $i = 1,...,n$, then
\begin{align}
(x_{0,1} - x_{n,n+1}, \cdots, x_{n-1,n} - x_{n,n+1}) = -(\delta_{1}, \cdots, \delta_{n})\Lambda_{n}. \label{for:x-x delta}
\end{align}
%Set element $e_{i} = x_{k(i),k(i)+1} - x_{i-sgnc_{i},i}$ for $i = 1,...,n$, which in fact is the $e_{i} $ in Fig. \ref{fig:The construction of surface F}(3). 
Notice that the relations of the columns of $\dot{W}$ are equivalent to
\begin{align}
    (t-1)e_{i} - sgnc_{i} \cdot t \cdot \delta_{i} = 0.\label{for:x-1 ei}
\end{align}
This implies that in the Alexander module, which is a torsion module, $(t-1)e_{i}$, $i = 1,...,n$ are equivalent to $\delta_{i}$, $i = 1,...,n$, which further implies that $e_{1}, \cdots, e_{n}$ form generators for the Alexander module.

We aim to show the following.

\textbf{Claim 3:} For the Alexander module, $R(t)$ is a presentation matrix, under generators $e_{1}, \cdots, e_{n}$. 

In fact, the relations of $\dot{W}$ are also equivalent to $e_{i} - t(x_{k(i),k(i)+1} - x_{i+sgnc_{i},i}) = 0$, that is,
\[
(e_{1},\cdots,e_{n})=t(x_{0,1}-x_{n,n+1},\cdots,x_{n-1,n}-x_{n,n+1})\dot{\Psi}^{*}.
\]
Combining (\ref{for:x-x delta}) and (\ref{for:x-1 ei}), we get
\begin{align}
 (e_{1},\cdots,e_{n})=-t(\delta_{1},\cdots,\delta_{n})\Lambda_{n}\dot{\Psi}^{*}=-(e_{1},\cdots,e_{n})(t-1)X\Lambda_{n}\dot{\Psi}^{*}.
\end{align}
We obtain a presentation matrix $I_{n}+(t-1)X\Lambda_{n}\dot{\Psi}^{*}$, which equals $-R(t)$ by (\ref{for:R(t) classicla knots}).
On the other hand, noticing that $\dot{\Psi}^{*}-\dot{\Phi}^{*}=\Lambda_{n}^{-1}X$, this matrix is $ -X\Lambda_{n}\dot{W}^{*}$.
\end{proof}

\subsubsection{Polyak-Viro formula}\label{sss:Polyak Viro formula}
As an application of Corollary \ref{cor:nomal Alexander polynomial classical knot}, we can derive Chmutov, Khoury and Rossi's formula\cite{ChmutovKhouryRossi2009} at least for the second coefficient $\mathbf{c}_2$ of Conway polynomial.

\begin{figure}[h]
\centering
\includegraphics[scale=0.3]{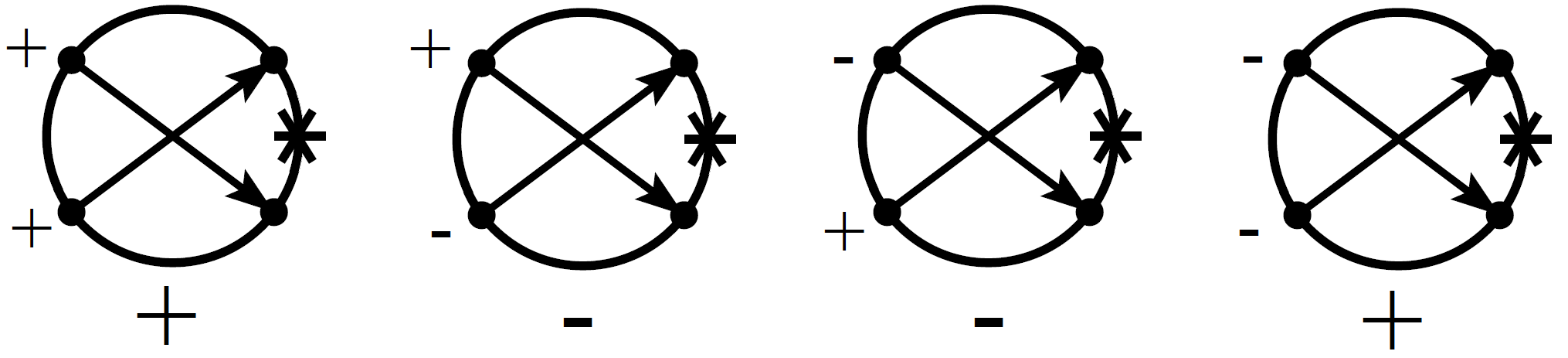}
\caption{Arrow diagrams and their signs to compute $\mathbf{c}_2$.}
\label{fig:arrow}
\end{figure}

\begin{cor}
Let $\mathbb{G}$ be a Gauss diagram for a knot $K$. Then $\mathbf{c}_2$ is the signed sum of the sub-Gauss diagrams of $\mathbb{G}$ as shown in Fig. \ref{fig:arrow}.
\end{cor}

\begin{rem}
The directions of the arrows in our formula is the opposite of those in the formula in \cite{ChmutovKhouryRossi2009}, with the same signs. It is clear, however, that the resulting $\mathbf{c}_2$ takes the same value.
\end{rem}

\begin{proof}
Let $\Delta(t)$ be the Conway-normalized Alexander polynomial for a knot, then we know that $\Delta(1)=1$, $\Delta'(1)=0$ and $\Delta''(1)=2\mathbf{c}_2$. So $\mathbf{c}_2=\frac{1}{2}|\mathcal{R}(t))''|\big |_{t=1}$. 

Let $X^{ij}$ be the matrix obtained by taking the derivatives of the $i$th and $j$th columns of $\mathcal{R}(t)$ and letting $t=1$. A simple calculation gives $\frac{1}{2}|\mathcal{R}(t))''|\big |_{t=1}=\Sigma_{i<j}|X^{ij}|$. As $\mathcal{R}(1)=I_n$, it is easy to see that $|X^{ij}|=|X^{ij}_{i,j;i,j}|$, where $X^{ij}_{i,j;i,j}$ is the $2 \times 2$ submatrix of $X^{ij}$ formed by the $i$th and $j$th rows and columns.

It remains to compute each $|X^{ij}_{i,j;i,j}|$. From (\ref{for:normal R(t) classical}) we can calculate
\begin{align}
    \mathcal{R}'_{ii}(1) &= 0.  \nonumber\\
    \mathcal{R}'_{ji}(1) &=\begin{cases}
-sgnP_i sgnc_j, & \text{if } c_{j0} \in P_i; \\
0, & \text{if } c_{j0} \notin P_i.
\end{cases} \quad \forall j \neq i. \nonumber
\end{align}

Using this formula, it is easy to complete the proof by discussing the value of $|X^{ij}_{i,j;i,j}|$ according to the directions and signs of the $i$th and $j$th arrows.
\end{proof}

\subsection{New formulae for virtual knots and long virtual knots.}

Let $\mathbb{K}$ be a diagram for a virtual knot $K$, with $n$ true crossings $c_{1},...,c_{n}$. Let $\mathbb{G}$ be the Gauss diagram of $\mathbb{K}$.

\subsubsection{New formulae for virtual and  long virtual knots }
Consider the matrix $R(t)$ defined by (\ref{for:R(t) classicla knots}). Note that Claim 1 in the proof of Proposition \ref{prop:relation with  Wirtinger} does not hold in general for virtual knots. This proof is just valid for virtual knots with a based point. So it is valid for \emph{long virtual knot}\cite{GoussarovPolyakViro2000}.
\begin{prop}
    For a long virtual knot, the Alexander polynomial is $ |R(t)|$, up to multiplication by $\pm t^{\pm k}$.
\end{prop}

For closed virtual knots, following the proof of Proposition \ref{prop:relation with  Wirtinger}, we can define a matrix
\begin{align}
 \tilde{R}(t) = - \left( 
 \begin{array}[]{cc}  
 X &    \\
   & 1
 \end{array}\right)
 \left( 
 \begin{array}[]{cc}  
 \Lambda_{n} &    \\
   & 1
 \end{array}\right)
\tilde{W} \nonumber
\end{align}
and obtain the following
\begin{prop}
    For a virtual knot, the Alexander polynomial is the greatest common divisor of the $(n - 1) \times (n - 1)$ minors of the matrix $\tilde{R}(t)$.
\end{prop}
However, this result looks not very interesting.

Nevertheless, we can say that, the Alexander polynomial of a virtual knot is a factor of the greatest common divisor of the determinants by taking base point in the interior of all the different over-passing arcs on $\mathbb{G}$.

\subsubsection{Contracted formula for long virtual knots.}
For virtual knot diagrams, there may be arrows with both contiguous starts and contiguous tips in the Gauss diagram $\mathbb{G}$. Merge all such adjacent arrows together and label each merged arrow $c_{i}$ with the number of the original arrows, denoted $\omega_{i}$, then we call it \emph{contracted Gauss diagram}, denoted $\mathbb{G}^{\omega}$. Given a based point $p$, we define a matrix $\boldsymbol{R}(t)$ from $\mathbb{G}^{\omega}$ as follows.
\begin{align}
\boldsymbol{R}_{ii}(t) &=\begin{cases}-1,&if\ \textit{sgn}P_{i}=\textit{sgn}c_{i};\\ -t,&if\ \textit{sgn}P_{i}=-\textit{sgn}c_{i}.\end{cases} \nonumber\\
\boldsymbol{R}_{ji}(t)&=\begin{cases}t^{\omega_{i}}-1,&if\ c_{j0}\in P_{i}\ and\ \textit{sgn}P_{i}=\textit{sgn}c_{j};\\ 1-t^{\omega_{i}},&if\ c_{j0}\in P_{i}\ and\ \textit{sgn}P_{i}=-\textit{sgn}c_{j};\\ \qquad\qquad 0,&if\ c_{j0}\notin P_{i}.\end{cases}\quad\forall j\neq i, \label{for:virtual contract}
\end{align}

An argument similar to the one used for ribbon graph shows that $|R(t)|=\pm|\boldsymbol{R}(t)|$. Thus we have 
\begin{prop}
    For a long virtual knot, the Alexander polynomial is $ |\boldsymbol{R}(t)| $, up to multiplication by $\pm t^{\pm k}$.
\end{prop}

\begin{cor}
The Alexander polynomial of a long virtual knot depends only on contracted Gauss diagram.
\end{cor}

This new formula (\ref{for:virtual contract}) can sometimes simplify the calculation of Alexander polynomial considerably.

\begin{figure}[h]
\centering
\includegraphics[width=0.5\textwidth]{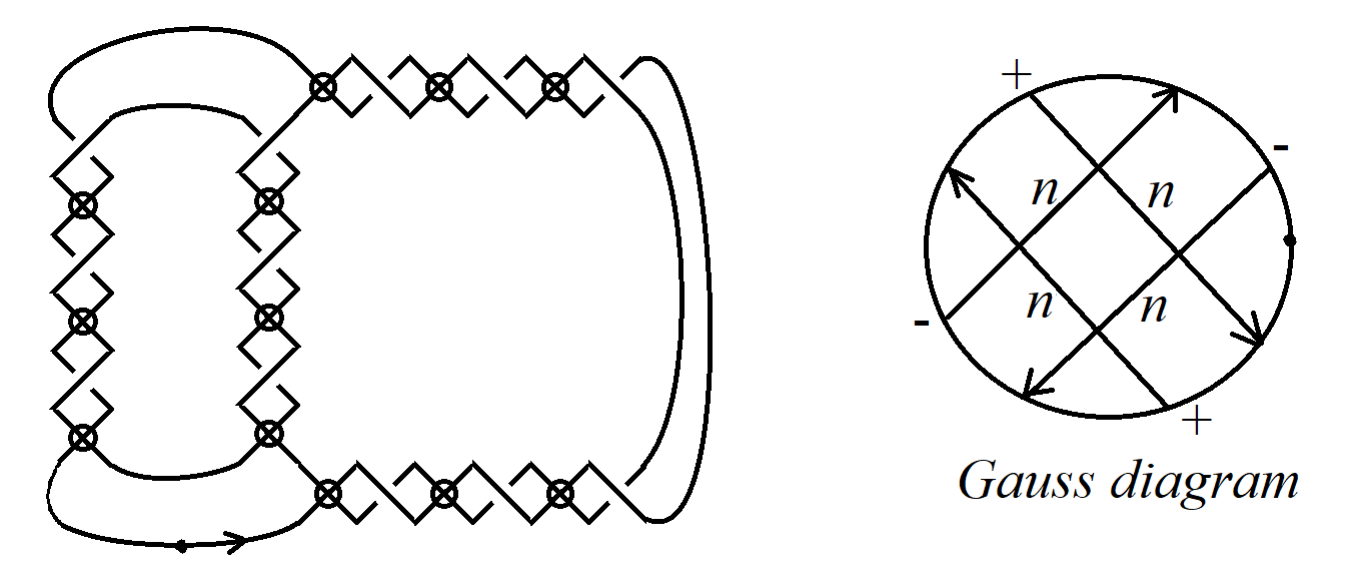}
\caption{A series of long virtual knots $K_{n}$. Here $n=3$.}\label{fig:long virtual knots}
\end{figure}

\begin{example}
For the long virtual knot $K_{n}$ in Fig. \ref{fig:long virtual knots}, we have
\begin{align}
    \boldsymbol{R}(t)=\left(\begin{array}[]{cccc}-t&0&0&0\\ 1-t^{n}&-1&1-t^{n}&0\\ t^{n}-1&1-t^{n}&-1&t^{n}-1\\ 0&t^{n}-1&0&-t\end{array}\right),
\end{align}
thus the Alexander polynomial $ \Delta_{K_{n}}(t)\dot{=} |\boldsymbol{R}(t)|\dot{=}-(t^{n}-1)^{3}+t^{n+1}(t^{n}-2)$. 
\end{example}

\section{Questions.}
We propose three questions, which may be not quite difficult, but of independent interests.

\subsection{On fusion number}

The number $g$ in the definition of ribbon diagram is called \textit{ribbon number} of the ribbon. However, the most researched ribbon knot invariant is \textit{fusion number}, which is the minimal number of fusions among all the ribbons for it. For a ribbon, fusion number can be much smaller than ribbon number. As ribbon graph is a reasonable notation to represent a ribbon, we ask

\begin{Question}
How to know the fusion number of a ribbon from its ribbon graph?
\end{Question}

\subsection{On half Alexander polynomials of symmetric unions.}

Symmetric unions were first introduced by Kinoshita and Terasaka \cite{KinoshitaTerasaka1957} and generalized by Lamm in \cite{Lamm2000}. Each symmetric union has a canonical ribbon $R$ by adding twists to the canonical ribbon of $\mathbb{K}\#-\overline{\mathbb{K}}$, where $\mathbb{K}$ is its partial knot. Lamm \cite{Lamm2000} conjectured in 2000 that every knot is a symmetric union, which is still open \cite{Lamm2012}.

For a symmetric union, the ribbon diagram of $R$ is the middle of Fig. \ref{fig:Gauss diagram and ribbon graph} with several small triangles reversed, which corresponds to changing several signs in the Gauss diagram of $\mathbb{K}$, and to changing the directions of several edges in the ribbon tree. Then we can use $W^{*}(t)$ in Theorem \ref{thm:algo3}, or equivalently $(n-1)\times(n-1)$ minor of $W$ from Fig. \ref{fig:FoxofWirtinger} with signs of some crossings changed, to compute $A_{R}(t)$.

For any knot, $(n-1)\times(n-1)$ minor of $W$ gives its Alexander polynomial, which is palindromic. We ask

\begin{Question}
For any $f(t)\in\mathbb{Z}[t]$ with $f(1)=\pm 1$, is there is a canonical ribbon of symmetric union so that the half Alexander polynomial $A_{R}(t)\doteq f(t)$?
\end{Question}

If the answer was false, one might easily get infinitely many counterexamples to Lamm's conjecture.

\subsection{On Polyak-Viro formula of the Conway polynomial}
Chmutov, Khoury and Rossi gave in \cite{ChmutovKhouryRossi2009} a Polyak-Viro formula for coefficients of the Conway polynomial. We have derived their formula for the second coefficients from our formula, in Subsection \ref{sss:Polyak Viro formula}. We ask

\begin{Question}
Can we derive all the Polyak-Viro formulas for coefficients of the Conway polynomial from our formula (\ref{for:normal R(t) classical}) ?
\end{Question}

Since the proof given by Chmutov, Khoury and Rossi was purely combinatorial, if so, we may obtain a topological understanding of their formulas.

\quad

\textbf{Acknowledgement:}
I am grateful to Jianhua Tu for his direction in graph theory and his help in writing Subsection 3.1 and 3.2. I gratefully acknowledges Professor Tetsuya Ito for his guidance in writing this paper. I thank Juhasz, Kauffman and Ogasa for their patience in discussing with me when I found there was a very implicit connection between our Theorem \ref{thm:classical} and their construction of twins Alexander polynomial in \cite{JuhaszKauffmanOgasa2022}.

\end{document}